\documentclass[11pt, leqno]{article}
\usepackage{graphicx}
\usepackage{amsfonts}
\usepackage{amsmath}
\usepackage{amsthm}
\usepackage{rotating}
\usepackage{verbatim}
\usepackage[export]{adjustbox}
\usepackage{oxford3}
\usepackage{varioref}
\usepackage{float}

\usepackage{color}
\usepackage{subcaption}
\usepackage{ amssymb }
\usepackage{hyperref}
\newcommand\vep{{\varepsilon}}

\newcommand\bA{{\bf A}}

\newcommand\mcD{{\mathcal D}}

\newcommand\bH{{\bf H}}
\newcommand\mcH{{\mathcal H}}

\newcommand\bI{{\bf I}}

\newcommand\bM{{\bf M}}

\newcommand\mcN{{\mathcal N}}

\newcommand\mcS{{\mathcal S}}
\newcommand\bs{{\bf s}}

\newcommand\bX{{\bf X}}

\usepackage{enumerate}

\newcommand\mbR{{\mathbb R}}

\newtheorem{corollary}{Corollary}
\newtheorem{lemma}{Lemma}
\newtheorem{theorem}{Theorem}
\newtheorem{definition}{Definition}                                                                  \newtheorem{assumption}{Assumption}
\theoremstyle{definition}

\newtheorem{examplex}{Example}
\newenvironment{example}{\pushQED{\qed}\examplex}{\popQED\endexamplex}

\begin{document}
\title{High-Dimensional Adaptive \\ Function-on-Scalar Regression }

\author{Zhaohu Fan\thanks{\small Departments of Industrial Engineering and Statistics, Penn State University,
University Park, PA, 16802} 
\qquad 
Matthew Reimherr\thanks{\small Department of Statistics,
Penn State University,
University Park, PA, 16802}\ \thanks{Corresponding Author: \href{mailto:mreimherr@psu.edu}{mreimherr@psu.edu}} \\ Pennsylvania State University
}

\date{}
\maketitle

\begin{abstract}
Applications of functional data with large numbers of predictors have grown precipitously in recent years, driven, in part, by rapid advances in genotyping technologies.  Given the large numbers of genetic mutations encountered in genetic association studies, statistical methods which more fully exploit the underlying structure of the data are imperative for maximizing statistical power.  However, there is currently very limited work in functional data with large numbers of predictors.  Tools are presented for simultaneous variable selection and parameter estimation in a functional linear model with a functional outcome and a large number of scalar predictors; the technique is called AFSL for \textit{Adaptive Function-on-Scalar Lasso}.  
It is demonstrated how techniques from convex analysis over Hilbert spaces can be used to establish a functional version of the oracle property for AFSL over any real separable Hilbert space, even when the number of predictors, $I$, is exponentially large compared to the sample size, $N$.  AFSL is illustrated via a simulation study and data from the Childhood Asthma Management Program, CAMP, selecting those genetic mutations which are important for lung growth.
\end{abstract}



\section{Introduction}

Many scientific areas are faced with the challenge of extracting information from increasingly large, complex, and highly structured data sets.  A great deal of modern statistical work focuses on developing tools for handling such data.  Networks, high dimensional data, images, functions, surfaces, or shapes, all present data structures which are not well handled under a traditional univariate or multivariate statistical paradigm.  
In this paper we present a new methodology, which we call \textit{Adaptive Function-on-Scalar Lasso}, AFSL, for analyzing highly complex functional outcomes alongside large numbers of scalar predictors.     
Such data is becoming increasingly common due to the prevalence of inexpensive genotyping technologies.  Genome-wide association studies, GWAS, examine hundreds of thousands or millions of genetic markers, attempting to find those mutations which are associated with some outcome or phenotype.  Many phenotypes of interest are now complex outcomes, such as longitudinal measurements or biomedical images. 

The functional linear model, FLM, is one of the primary modeling tools in FDA.  There, one assumes that the outcome is linearly related to some set of predictors.
FLM are often categorized by whether the outcome, the predictor, or both is functional \cite{reiss:etal:2010}.  While the literature on FLM is now vast, \citetext{morris:2015} outlines most of the key work on low dimensional FLM.  However, to date, relatively little has been done when one has a large number of predictors relative to the sample size, and of the work that does exist, nearly all of it is for scalar-on-function FLM, the opposite setting we consider \cite{MaKo:2011,Lian:2013,gertheiss:maity:staicu:2013,fan:james:rad:2015}. 
For the function-on-scalar setting, we are aware of only two other works.  \citetext{chen:goldsmith:ogden:2016} consider a basis expansion approach with an MCP style penalty and fixed number of covariates.  They also use a \textit{pre-whitening} technique to exploit the within function dependence of the outcomes.  Asymptotic theory is developed for a fixed number of predictors and basis functions.
\citetext{M1} developed the Function-on-Scalar LASSO, FSL, which allows for an exponential number of predictors relative to the sample size and establishes optimal convergence rates of the estimates.  In particular, it was shown that FSL achieves the same rates of convergence as in the scalar case.  However, as in the scalar case, this approach leads to estimates with a non--negligible asymptotic bias due to the nature of the penalty.  

To alleviate the bias problem inherent in FSL, we propose here an adaptive version, AFSL.  In addition to providing a novel statistical methodology, we develop a new theoretical framework needed to establish its asymptotic properties.   In particular, functional subgradients and tools from convex analysis over Hilbert spaces are needed.  In contrast, theory for FSL is built entirely on \textit{basic} inequalities and \textit{concentration} inequalities, no theory for subgradients was required.  The contributions of this paper are thus as follows (1) we define a new variable selection and estimation tool, AFSL, which alleviates the bias problems of FSL (2) we demonstrate how several tools and techniques from convex analysis can be used for functional data problems and (3) we define a functional version of the oracle property and show that AFSL achieves it.   Additionally, we also go a step beyond the traditional oracle property which states that the estimates recover the correct support and are asymptotically normal, by showing that the oracle estimate and the AFSL estimate are actually asymptotically equivalent.

The remainder of the paper is organized as follows.  In Section 2, we provide the necessary background material. In section 3 we outline the AFSL framework and in section 4 we establish the oracle property.  A simulation study and an application to the Childhood Asthma Management Program, CAMP, is given in Section 5.  Concluding remarks are given in Section 6.  All theory is provided in the Appendix.

\section{Background}\label{s:back}
Let $\mcH$ denote a real separable Hilbert space, $\langle \cdot, \cdot \rangle$ its inner product, and $\| \cdot \|$ the inner product norm.  In a function-on-scalar linear model we have that
\begin{equation}
Y_{n}=\sum_{i=1}^{I} X_{ni}\beta_{i}^\star+\vep_{n} = X_n^\top \beta^\star + \vep_n \quad
1 \leq n \leq N, 
\end{equation} 
where $Y_{n} \in \mcH$ is a functional outcome, $\beta_i^\star \in \mcH$ is a functional regression parameter, $\vep_n \in \mcH$ is a functional error, and $X_{ni} \in \mbR$ is a scalar predictor.  Throughout, we will use a $\star$ to denote the true data generating parameter so as to distinguish from $\beta$, which will usually represent a dummy variable or the argument of a function.  The most commonly encountered space for $\mcH$ is $L^2[0,1]$, i.e., the outcome is a function of time, though other spaces are used as well including spatial domains or product spaces (for functional panels or multivariate functional data).  If one wants to incorporate smoothness assumptions of the data then one can work with Sobolev spaces or Reproducing Kernel Hilbert Spaces, RKHS.  We provide a few examples here to help emphasize the functional nature of the data and highlight the wide impact of our theory and methods.
\begin{example}
Let $\mcH = L^2(\mcD)$, where $\mcD$ is a compact subset of $\mbR^d$.  Then we write the model
\[
Y_{n}(\bs)=\sum_{i=1}^{I} X_{ni}\beta_{i}^\star(\bs)+\vep_{n}(\bs)
\qquad \bs \in \mcD,
\]
and the norm is written as
\[
\| x\|^2 = \int_\mcD x(\bs)^2 \ d \bs.
\]
\end{example}
\begin{example}
Let $\mcH = H^{1,2}[0,1]$, i.e., the Sobolev space of real valued functions over the unit interval with one square integrable derivative.
 Then we have
\[
Y_{n}(t)=\sum_{i=1}^{I} X_{ni}\beta_{i}^\star(t)+\vep_{n}(t) 
\qquad t \in [0,1], 
\]
with the norm given by
\[
\| x\|^2 = \int_0^1 x(t)^2 \ dt +  \int_0^1 x^\prime (t)^2 \ dt.
\]
\end{example}

\begin{example}
Let $\mcH$, be an RKHS of functions over the unit interval with kernel operator $K$.
 Then we have
\[
Y_{n}(t)=\sum_{i=1}^{I} X_{ni}\beta_{i}^\star(t)+\vep_{n}(t), 
\]
but now the norm becomes
\[
\| x\|^2 = \langle K^{-1} x , x \rangle,
\]
where $\langle \cdot, \cdot \rangle$ is the $L^2[0,1]$ inner product and $K^{-1}$ is the inverse (operator) of $K$.
\end{example}

The most common method for estimating $\beta$ is least squares, i.e. minimizing
\[
L_0(\beta) = \frac{1}{2}\sum_{n=1}^N \| Y_n - X_n^\top \beta \|^2.
\]
However, when $ I > N$ the problem is ill-posed and a unique minimizer does not exist.  To address this issue, \citetext{M1} introduced FSL, which is defined as the minimizer of
\[
L_\lambda(\beta) = \frac{1}{2} \sum_{n=1}^N \| Y_n - X_n^\top \beta \|^2 + \lambda \sum_{i=1}^I \| \beta_i\|.
\]
The term $\sum_{i=1}^I \| \beta_i\|$ can be viewed as a type of $\ell_1$ norm on the product space $\mcH^{I}$ which induces a sparse estimate.  As the name implies, this approach is an extension of the classic Lasso \cite{tibshirani:1996} to functional outcomes.  It was shown that under appropriate conditions, the FSL estimate achieves optimal convergence rates and that these rates are the same regardless of the Hilbert space.  In other words, the converge rate for $Y \in \mbR$ is the same as for $Y \in L^2[0,1]$.

As in the scalar case, FSL does not achieve the oracle property.  This is due to the non-negligible bias of the estimates.  In the next section we define an adaptive version of FSL and show that it achieves a type of oracle property.  Interestingly, we will use a different formulation of the oracle property than the one commonly used in the scalar setting, namely, that estimates are asymptotically normal.  In the functional setting, some estimates are not asymptotically normal even in the low dimensional setting \cite{cardot:2007}, and thus removing this connection is useful.
  
\section{Adaptive Function-on-Scalar Regression}
We now define an adaptive version of FSL we call AFSL.  As with FSL, this can be viewed as an  extension of the adaptive lasso  \cite{zou:2006} to functional outcomes.  To produce an adaptive estimate, we introduce data driven weights into the target function.  Namely, we denote by $\hat \beta$ the AFSL estimate, which is the minimizer of
\begin{align}
L_\lambda^A(\beta) = 
\frac{1}{2}\sum_{n=1}^N \| Y_n - X_n^\top \beta\|^2 + \lambda \sum_{i=1}^I \tilde w_i \| \beta_i\|. \label{e:afsl}
\end{align}
These weights can be chosen in a variety of ways, but the approach we take here is to run FSL to obtain preliminary estimates $\tilde \beta_i$.  We then set $\tilde w_i = \| \tilde \beta_i\|^{-1}$ so that small estimates are shrunk to zero and large estimates are not shrunk as substantially.  It is possible that $\tilde w_i = \infty$ in which case the predictor is dropped.  This has the advantage that running AFSL after FSL is usually very fast as most of the variables have already been screened.  \textcolor{black}{Another option, which is used in the scalar case by \cite{h08a}, is to do marginal ordinary least squares regression to obtain each weight.  However, this does not reduce the dimension of the problem in anyway and requires an additional theoretical justification, thus we do not pursue it further here.  }

\citetext{M1} developed optimal asymptotic convergence rates of FSL by combining a basic inequality with concentration inequalities and a restricted eigenvalue condition.  Here we aim to show that AFSL achieves the oracle property which is a much stronger property.  It is therefore useful to switch to a different technique which involves functional subgradients.  
Given that many readers may not be familiar with such concepts for general Hilbert spaces, we provide a brief reminder of the core concepts.  We refer interested readers to \citetext{boyd:2004,bauschke:2011,shor:2012} for more details on subgradients and convex analysis, and \citetext{barbu:precupanu:2012} specifically for convex analysis over Hilbert spaces.
\begin{definition} \label{d:sub}
Let $f:\mcH \to \mbR$ be a convex functional. An element $h \in \mcH$ is called a subgradient of $f$ at $x$ if for all $y \in \mcH$ we have
\begin{equation*}
f(y)-f(x)\geq \langle h,y-x\rangle.
\end{equation*}
\end{definition}
We emphasize that in definition \ref{d:sub} the points $x$, $y$, and $h$ are elements of an arbitrary Hilbert space.  If $\mcH = \mbR$, then they are scalars, if $\mcH = \mbR^d$ then they are vectors, and if $\mcH = L^2[0,1]$ or any of the function spaces from \ref{s:back}, then they are functions.  Note that subgradients are unique if the functional is differentiable, but they also exist for nondifferentiable functions in which case they need not be unique.  
\begin{definition}
The set of subgradients of $f$ at $x$ is called the subdifferential and denoted $\partial f(x)$
\end{definition}

Thus, one has a generalized notion of derivative for convex functionals.  One can immediately see from Defintion \ref{d:sub} that if $0$ is in the subdifferential at a point $x$, then that point is a minimizer.   Both the FSL and AFSL are convex, we thus have the following property.

\begin{theorem} \label{t:subgrad} The subdifferential of $L^{A}_{\lambda}\left(\beta\right)$ with respect to $\beta_i$ is given by
\[
\frac{\partial}{\partial \beta_{i}}  L_{\lambda}^{A}(\beta_{i})=-\sum_{n=1}^{N}X_{ni}(Y_{n}-X_{n}^{T}\beta)+\lambda \tilde w_{i}
\begin{cases}
\beta_{i}\|  \beta_{i}\|^{-1} & \text{if } \| \beta_{i}\|\neq 0 \\
\{h\in \mcH:~~ \|h\|\leq 1\}& \text{if } \|  \beta_{i}\| =  0
\end{cases}.
\]
\end{theorem}
We again emphasize that $\beta_i$ is an element of $\mcH$ (i.e.\ is a function), thus the subdifferential in Theorem \ref{t:subgrad} is a set of objects from $\mcH$.
Using Theorem \ref{t:subgrad} we can't, in general, obtain an explicit expression for $\hat \beta$, but we can characterize the solution enough to establish the oracle property.  In particular, $\hat \beta$ must satisfy the following corollary.
\begin{corollary}\label{c:sub}
If $\hat \beta$ satisfies
\[
-\sum_{n=1}^{N}X_{ni}(Y_{n}-X_{n}^{\top}\hat \beta)+\lambda \tilde w_{i} \hat \beta_{i}\|  \hat \beta_{i}\|^{-1} = 0, \qquad \text{for $\hat \beta_{i} \neq 0$}
\]
 and
\[
\left \| \sum_{n=1}^{N}X_{ni}(Y_{n}-X_{n}^{\top}\hat \beta) \right \|
< \lambda \tilde w_i
\qquad \text{for $\hat \beta_{i}= 0$},
\]
then $\hat \beta$ is a minimizer of $L_A(\beta)$.
\end{corollary}
Using these properties combined with appropriate assumptions on the errors and the design matrix, we establish the oracle property in the next section.  However, to help better understand the structure of AFSL, we consider the orthogonal design case as an illustrative example.  In that case, both FSL and AFSL produce explicit estimates.

\begin{example}[Orthogonal Design]
 For this example only assume that $I < N$ is fixed and that $N^{-1} \mathbf{X^\top} \mathbf{X} = \mathbf{I}_{I \times I}$. 
 From Corollary \ref{c:sub} if $\hat \beta_i$ is not zero then it must satisfy
\[
\sum_{n=1}^N X_{in} Y_n -   N\hat \beta_i = \lambda \tilde w_i \hat \beta_i \| \hat \beta_i\|^{-1}.
\]
Since the design is orthogonal, $\hat \beta_i^{LS} =  \frac{1}{N}\sum_{n=1}^N X_{in} Y_n $ is the least squares estimator.  
So we have
\[
 \hat \beta_i = \frac{1}{1+\lambda N^{-1} \tilde w_i \| \hat \beta_i\|^{-1}} \hat \beta_i^{LS}.
\]
We can get the norm of $\hat \beta_i$ by considering
\[
 \|\hat \beta_i^{LS} \| = (1+\lambda N^{-1}  \tilde w_i  \| \hat \beta_i\|^{-1} ) \| \hat \beta_i\|
=  \| \hat \beta_i\| +\lambda N^{-1} \tilde w_i,
\]
which implies that
\[
\| \hat \beta_i\| = \| \hat \beta_i^{LS}\|     - \lambda N^{-1} \tilde w_i.
\]
After a little algebra, the AFSL estimate can be expressed as
\begin{align*}
\hat \beta_i = \left(1   -  \frac{\lambda \tilde w_i }{ N \| \hat \beta_i^{LS} \|} \right)\hat \beta_i^{LS}.
\end{align*}
Turning to the case where $\hat \beta_i = 0$ we have that
\[
\|  X_i^\top (Y - \bX \hat \beta) \| < \lambda \tilde w_i,
\]
or equivalently
\[
\| \hat \beta_i^{LS} \| < {\lambda N^{-1} \tilde w_i}.
\]
So, as in the scalar case, the FSL and AFSL procedures can be viewed as a soft thresholding applied to the least squares estimator:
\[
\hat \beta_i
=  \left(1   -  \frac{\lambda \tilde w_i }{ N \| \hat \beta_i^{LS}\|} \right)^+ \hat \beta_i^{LS}. 
\]  
At this point, we can now establish conditions for the oracle property to hold.  Notice that for $\beta_i = 0$, we want  $\lambda \tilde w_i N^{-1}   \| \tilde\beta_i\|^{-1}$ to grow quickly, while for $\beta_i \neq 0$, we want it to shrink quickly (so as to reduce the bias as much as possible). If $\beta_i =0$ then $\hat \beta_i^{LS} = O_P(N^{-1/2})$ and we can assume the same for $\tilde w_i^{-1}$ (e.g. this would be true if we used the least squares estimates to compute $\tilde w_i$).  We have
\[
\frac{\lambda \tilde w_i }{  N \| \hat \beta_i^{LS}\|} ={\lambda} O_P(1),
\]
so as long as $\lambda  \to \infty$, then $\hat \beta_i$ will be shrunk to zero.  Conversely, if $\beta_i \neq 0$, then 
\[
\frac{\lambda \hat w_i }{  N \| \hat \beta_i^{LS}\|} = {\lambda}N^{-1} O_P(1),
\]
so as long as $\lambda/N \to 0$, then the estimate will be asymptotically unbiased.  Under these conditions, AFSL has variable selection consistency, but we need a bit more control of $\lambda$ to ensure that the estimate is asymptotically equivalent to the oracle.
When $\beta_i \neq 0$ we have
\[
N^{1/2} \| \hat \beta_i -\hat \beta_i^{LS}\| = 
N^{-1/2} \lambda O_P(1).
\]
So, if $\lambda N^{-1/2} \to 0$ then the LS estimator for the nonzero coordinates (this is usually called the oracle estimator) is asymptotically equivalent to the AFSL estimator and thus is asymptotically normal.   This is because, for $\beta_i \neq 0$
\[
N^{1/2}(\hat \beta_i - \beta_i^\star ) =  N^{1/2}(\hat \beta_i^{LS} - \beta_i^\star ) +  N^{1/2}(\hat \beta_i - \hat \beta_i^{LS} )
=  N^{1/2}(\hat \beta_i^{LS} - \beta_i^\star ) + o_P(1).
\]
Therefore, we can conclude that for fixed $I$ and an orthogonal design, AFSL achieves the oracle property as long as $1 \ll \lambda \ll  N^{1/2}$.  Later on, we will see that this main dynamic still holds.  With a ``looser" control of $\lambda$, we can get selection consistency, but to ensure AFSL is equivalent to the oracle estimate, we need slightly tighter control.

\end{example}

\section{Oracle Property}
We begin by making the FLM assumption more explicit.
\begin{assumption} \label{a:model}
Let $Y_{1},...,Y_{N}$ be independent random elements of a real separable Hilbert space, $\mathcal{H}$, satisfying the functional linear model:
\begin{equation}
Y_{n}=\sum_{i=1}^{I}X_{ni}\beta_{i}^\star+\epsilon_{n}.
\end{equation}
Assume the N$\times$I design matrix $\mathbf{X}=\{X_{ni}\}$ is deterministic and has standardized columns, and that $\epsilon_{n}$ are i.i.d. Gaussian random elements of $\mathcal{H}$ with mean function $0$ and covariance operator $C$. 
\end{assumption}
\noindent As is common in high dimensional settings, the Gaussian assumption is not crucial and our arguments can be readily generalized to errors with thicker tails.  We utilize the Gaussian assumption here to make the rates in our subsequent assumptions simpler and more interpretable.  We emphasize that the i.i.d.\ assumption means that response curves from different subjects are independent, but observations from the same subject are allowed to be dependent, and no restriction is placed on this dependence.  

Define the true support as $\mcS = \{ i \in {1,\dots, I} : \beta_{i}^\star\neq 0\}$ and $I_0 =|\mcS|$, i.e.\ the cardinality of $\mcS$.  
For notational simplicity, we will assume that $\bX$ is ordered such that it can be partitioned into
\[
\bX = ( \bX_1 \quad \bX_2),
\]
where $\bX_1 \in \mbR^{N \times I_0}$ are covariates with nonzero coefficients and $\bX_2 \in \mbR^{N \times (I - I_0)}$ are the covariates with zero coefficients.  

Next we define a functional version of the oracle property.  For scalar outcomes, this is usually phrased as having asymptotically normal estimates and estimating the true support with probability tending to one \cite{fan:li:2001}.  For functional data, we will divorce the normality from the oracle property because there are examples of functional estimates which are not asymptotically normal, e.g. scalar-on-function regression \cite{cardot:2007}.  Instead, we make the more direct assumption that the oracle property means that an estimator asymptotically has the correct support, and is asymptotically equivalent to the oracle estimator.  The oracle estimate in this case is defined as $\hat \beta_O = \{(\bX_{1} \bX^\top_{1})^{-1} \bX^\top _1 Y, {\bf 0}\}$, that is, all estimates from the second group are set to zero and not used to form estimates for the first group.
We now introduce the functional oracle property.
\begin{definition}  We say that an estimate $\hat \beta$ of $\beta$ satisfies the functional oracle property if
\begin{enumerate}
\item $P(\hat \beta \overset{s}{=} \beta) \to 1$ and
\item $\sqrt N (\hat \beta - \hat \beta_O) = o_P(1)$.
\end{enumerate}
\end{definition}
\noindent Here $\overset{s}{=}$ means that the two have the same support, i.e., same nonzero coordinates.  We now introduce the remaining more technical assumptions.  All of these conditions are commonly found in the literature and relate the orders of various terms.  
%
%
\begin{assumption}\label{a:tech}
We assume that following four conditions hold:
\begin{enumerate}
\item {\bf Minimum Signal:} Let $b_N = \min_{i=1,\dots,I_0} \| \beta_i\|$.  Then we assume that
$$ b_N^2 \gg \frac{I_0^2 \log (I) }{N}.$$
\item {\bf Tuning Parameter:} We assume the tuning parameter satisfies
$$ I_0^{1/2} \log(I) \ll \lambda \ll \frac{N b_N^2}{I_0^{1/2}}.$$ 
\item {\bf Design Matrix:} Let $\hat \Sigma_{1} = N^{-1} \bX_1^\top \bX_1$ , then we assume that the smallest and largest eigenvalues of $\hat \Sigma_1$ satisfy ,
\begin{equation*}
\frac{1}{\tau_{1}}\leq \sigma_{\text{min}}\left( \hat \Sigma_{1}\right)\leq \sigma_{\text{max}}\left(\hat \Sigma_{1}\right)\leq\tau_{1}
\end{equation*}
where $\tau_1 >0$ is fixed scalar.
\item {\bf Irrepresentable Condition:} Let $\hat \Sigma_{21} = N^{-1} \bX_2^\top \bX_1$, then we assume that
$$\|\hat \Sigma_{21}\hat \Sigma_{11}^{-1}\|_{op} \leq \phi<1,$$
where $\| \cdot \|_{op}$ is the operator norm and $\phi$ is a fixed scalar.
\end{enumerate}
\end{assumption}
Each of the above assumptions is common in the high dimensional regression literature.  The minimum signal condition allows the smallest $\beta_i$ to vary with the sample size, but it cannot get too small.  The tuning parameter assumption states a familiar rate, namely, the tuning parameter must grow, but not too quickly.  The eigenvalue assumption requires that the true predictors are well behaved and not extremely highly correlated.  The last assumption, the Irrepresentable condition, essentially says that the true and null predictors cannot be too correlated.  For the scalar oracle property it has been shown that this is essentially necessary and cannot be weakened \cite{zhao:yu:2006}.   	

The above assumptions, as we will see, will guarantee that the AFSL estimate is variable selection consistent.   However, to control the asymptotic bias, and thus ensure the AFSL estimate is asymptotically equivalent to the oracle, we also need the following.
\begin{assumption} \label{a:tech2}
The tuning parameter satisfies
\[
\lambda \ll \frac{\sqrt{N} b_N}{I_0^{1/2} }.
\]
\end{assumption}	
When we divide the right hand side of Assumption \ref{a:tech}.2 by the left hand side , and assuming that $b_N$ is bounded, we see the familiar assumption that
\[
\frac{N}{I_0 \log (I) } \to \infty.
\]
However, to control the bias enough for the full oracle property to hold, we also need Assumption \ref{a:tech2}.  If we take the right hand side of Assumption \ref{a:tech2} and divide by the left hand side of Assumption \ref{a:tech}.2 we arrive at another familiar form
\[
\frac{\sqrt{N}}{I_0 \log (I)} \to \infty.
\]
These rates are well known and necessary for scalar linear models \cite{buhlmann:2011,cai:2015}.  What is especially interesting here, is that we obtain the exact same rates for any separable Hilbert space.  In other words, the outcome can come from any separable Hilbert space, and the oracle property holds with exactly the same rates as in the scalar case. 	
	
We are now in a position to state our main result.	
\begin{theorem}(Functional Oracle property) \label{t:oracle}
Let $\tilde \beta$ be the FSL estimator and let $\hat{\beta}$ be a minimizer of \eqref{e:afsl} with weights $\tilde w_i = \|\tilde \beta_i\|^{-1}$.  We then have the following.  \begin{enumerate}
\item If Assumptions \ref{a:model} and \ref{a:tech} hold then
$$
P(\hat \beta \overset{s}{=} \beta) \to 1.
$$
\item  If in addition Assumption \ref{a:tech2} holds then
$$\sqrt{n}(\hat \beta  - \tilde\beta_{O}) = o_P(1).$$
\end{enumerate}
\end{theorem}
As a corollary, since the oracle estimate for function-on-scalar regression is asymptotically normal, we have the following.
\begin{corollary}
If Theorem \ref{t:oracle} holds, then for any nonzero sequence $x_N \in \mbR^{I_0}$ we have that
\[
\frac{1}{\sqrt{x_N^\top \hat \Sigma_1^{-1} x_N} }\sum_{i=1}^{I_0} x_{Ni} (\hat \beta_i - \beta_i^\star)
\overset{\mcD}{\to} \mcN(0, C),
\]
where $\mcN(0,C)$ is an $\mcH$ valued Gaussian process with zero mean and covariance operator $C$.
\end{corollary}

\section{Numerical Studies}
In this section we perform a simulation study to illustrate the performance gain of AFSL over FSL.  We also apply our methodology to data from Childhood Asthma Management project, a longitudinal genetic association study.  \textcolor{black}{However, before exploring these examples, we discuss how the methodology is implemented in both cases.  The data are first preprocessed into functional units via penalized (on the second derivative) bsplines basis expansions using the {\tt FDA} package in {\tt matlab}.  We use 100 cubic bsplines with equally spaced knots, and the smoothing parameter is chosen via generalized cross-validation.  
This processing is now standard in the FDA literature and we reference \citetext{h12} for more details.  After the functional units are constructed, we rotate the data to the functional principal component, FPC, basis.  This serves the dual purpose of working with an orthonormal basis which is also tailored to the data.  Practically, one can work with a relatively small FPC basis compared to bsplines.  We emphasize that we do not use this basis for the purposes of serious dimension reduction, so we take the number of FPCs to explain $\geq 99\%$ of the variance so that the FPCA and bsplines approximations are nearly identical.  Next, we rephrase the problem as a group lasso problem, as in \citetext{M1}, and use alternating direction method of multipliers, ADMM, to find the solution \cite{boyd:2011}.  This approach can be used for both FSL and AFSL by using the weights to scale the predictors instead of the $\beta$; one performs a change of variables with $\alpha_i = \tilde w_i \beta_i$, uses the discussed method for finding the $\hat \alpha_i$, and then changes variables back to the $\hat \beta_i$.  We choose the tuning parameter via BIC, though in the data application we also compare to the extended BIC as was used in \citetext{M1}.  A complete {\tt matlab} function and example are available for download from the corresponding author's website,
\url{www.personal.psu.edu/~mlr36}, 
which can be readily used by other researchers.  Access to the discussed CAMP data is free, but one must submit an application through NIH's dbGaP.  Those interested in the {\tt matlab} code applied to the CAMP data should contact the corresponding author directly.}

\subsection{Empirical Study}\label{s:sims}
We generate random samples of size $N = 100$ and $500$.  For each sample, we generate the nonzero $\beta_i$ from from a mean zero Gaussian process with covariance given in \eqref{e:fff1}.  Note that this is a Mat\'ern process with smoothness parameter $\nu = 2.5$, which means the $\beta$ are twice differentiable.  We also take the point-wise variance to be $\sigma^2 = 1$ and the range (i.e.\ how quickly the within curve dependence decays) to be $\tau_1 =0.25$.  
The errors, $\vep_n$, are generated in nearly the same way, but we change $\nu$ to $1.5$, meaning that the errors are rougher than betas (as would be natural in practice) possessing only one derivative.  The corresponding covariance function is given in \eqref{e:fff2}.  We also consider several different range parameters, $\tau_2=0.01,0.25, 1, 10$, to examine the effect of the within subject correlation. 

\textcolor{black}{
\begin{equation}
C(t,s) = \sigma^2 \Bigg(1 + \frac{ \sqrt{5}d }{\tau_1} +\frac{ 5d^2}{3 \tau_1^2 }   \Bigg) \exp \Bigg(-\frac{\sqrt{5}d}{\tau_1} \Bigg),
\quad d = |t-s|, \quad 
\nu= \tfrac{5}{2}.\label{e:fff1}
\end{equation}
}

\textcolor{black}{
\begin{equation}
C(t,s) = \sigma^2 \Bigg(1 + \frac{ \sqrt{3}d }{\tau_2} \Bigg) \exp \Bigg(-\frac{\sqrt{3}d}{\tau_2} \Bigg), \quad 
d = |t-s|,\quad 
\nu= \tfrac{3}{2}. \label{e:fff2}
\end {equation}}

\begin{figure}
\centering
\includegraphics[scale=0.35]{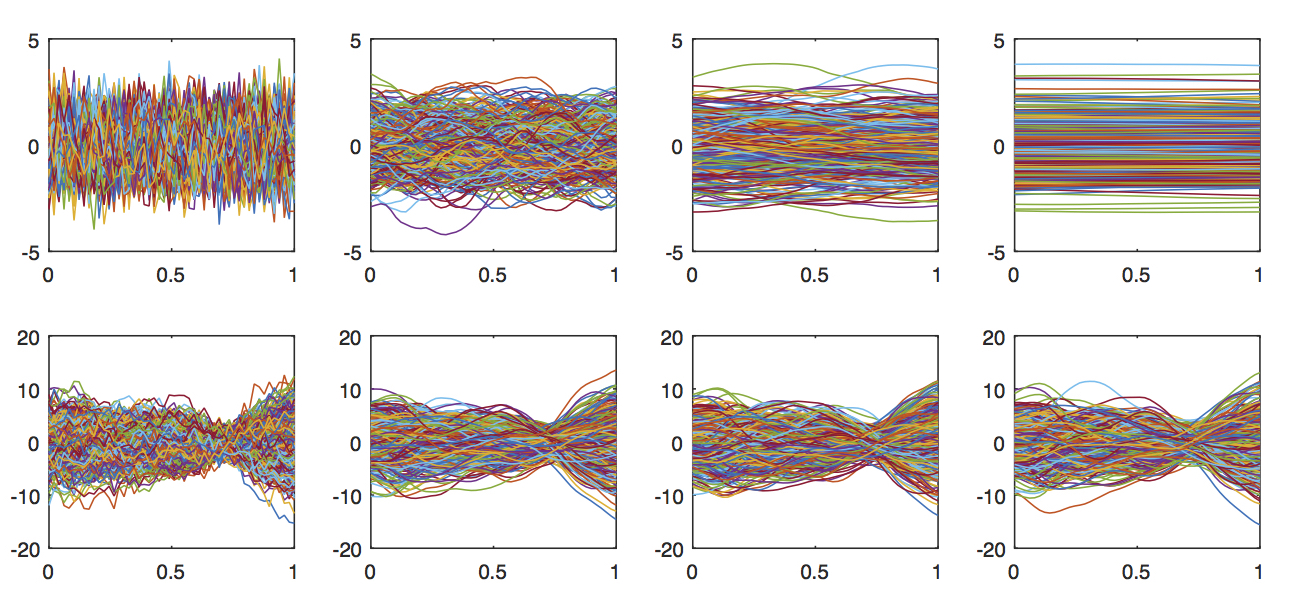}
\caption{Plots of error terms, from left to right, with range parameters $\tau_2$=0.01, 0.25, 1 and 10, on the top row and the corresponding outcome $Y$ on the bottom row.}
\label{tab:F111}
\end{figure}

Each curve is sampled at 50 evenly space points between $0$ and $1$.  The total number of predictors $I$ is equal to 500 while the total number of true predictors, $I_0$, is set to 10.  The predictors are $X_{ni}$ are taken to be standard normal with a time series style correlation structure: $Cov(X_{ni},X_{nj}) = \rho^{|i - j|}$. \textcolor{black}{As an illustration, in Figure \ref{tab:F111} we plot four examples of error terms, $\vep_n$, with four different range parameters ($\tau_2$=0.01, 0.25, 1 and 10) and the corresponding outcomes, $Y_n$, when $\rho = 0$}. 
We consider three different values of $\rho$ to examine the effect of predictor correlation on performance: $\rho = 0, 0.5, 0.99$ corresponding to low, medium, and high correlation.  We consider 1000 repetitions of each scenario and in every setting, we choose the tuning parameter $\lambda$ using BIC as in \citetext{M1}.
The results are summarized in Table \ref{tab:111}, where we report the average number of true positives, average number of false positives, average root-mean-squared prediction error, and the average computation time.  Each of these values is reported for both FSL and AFSL.

\begin{table}[ht]
\centering
\resizebox{5.8in}{!}{%
\begin{tabular}{ccccccccl} \hline
& \multicolumn{2}{c}{True Positives} &  \multicolumn{2}{c}{False Positives}& \multicolumn{2}{c}{RMSP} &  \multicolumn{2}{c}{ Time (seconds)}\\ \hline
      Parameters & {FSL} & {AFSL}& FSL &  AFSL& {FSL} & {AFSL}& FSL &  AFSL\\ \hline
       $\tau=0.01; \rho$=0; N=100 &   9.82	&9.78&	3.34&	0.42 &  0.0129&	0.0113	&156.9379	&165.5910   \\ 
     $\tau=0.01;\rho$=0; N=500&  9.99&	9.98&	17.92&	0.06&0.0051	&0.0050&	283.1586	&311.3364\\ 
      $\tau=0.01;\rho$=0.5; N=100   & 9.91&	9.84&	3.35	&0.41&0.0126&	0.0113	&172.7202&	180.0085\\
       $\tau=0.01;\rho$=0.5; N=500 &9.99	&9.96	&13.16	&0.03&   0.0049	&0.0051&	284.8761	&316.0269\\
             $\tau=0.01;\rho$=0.99; N=100 & 5.66&	5.18	&10.63	&3.36 & 0.0182&	0.0105&	965.0788	&859.6725\\ 
       $\tau=0.01;\rho$=0.99; N=500 & 9.82	&9.36&	16.68&	0.02 & 0.0037	&0.0034&	1167.7435	&1191.7981\\  \hline
$\tau=0.25; \rho$=0; N=100 & 8.62	&8.43	&1.65	&0.68 & 0.0193	&0.0128&	88.5702	&91.2917 \\  
  $\tau=0.25;\rho$=0; N=500&9.99	&9.87	&2.82&	0.13&    0.0054&	0.0051	&341.4661	&369.8791\\
      $\tau=0.25;\rho$=0.5; N=100   &8.25	&8.18&	1.82&	0.71& 0.0202	&0.0125	&98.3453	&98.1497\\
       $\tau=0.25;\rho$=0.5; N=500 & 9.99	&9.89&	2.52	&0.16&  0.0054&	0.0051	&367.4788	&395.9227\\
             $\tau=0.25;\rho$=0.99; N=100 &1.18	&1.05	&3.27	&2.14&0.0328	&0.0092	&996.5986	&476.5383 \\
       $\tau=0.25;\rho$=0.99; N=500 &8.82&	8.08&	15.98	&2.82 &0.0045&	0.0040&	1291.9321&	1310.9210\\ \hline
$\tau=1; \rho$=0; N=100 &  7.02	&6.91	&1.13	&0.52  &  0.0239	&0.0127&	65.5368&	62.4991  \\ 
  $\tau=1;\rho$=0; N=500& 9.96&	9.86	&2.09	&0.47&    0.0056&	0.0051&	305.7731	&332.5053\\
      $\tau=1;\rho$=0.5; N=100   & 6.34	&6.28&	1.74	&0.46&   0.0262&	0.0113&	70.4651&	58.5442 \\
       $\tau=1;\rho$=0.5; N=500 &9.99	&9.88	&2.51	&0.62&   0.0057	&0.0053&	308.2203&	333.3216\\
             $\tau=1;\rho$=0.99; N=100 &  6.08&	6.02	&1.06	&0.58& 0.0263&	0.0120	&62.9913	&54.8344\\
       $\tau=1;\rho$=0.99; N=500 &  9.98	&9.94	&1.96	&0.34 &0.0057&	0.0052&	297.8822	&322.7747 \\ \hline
$\tau=10; \rho$=0; N=100 &   7.02	&6.93	&1.31	&0.74 & 0.0239	&0.0123	&57.8824	&55.3399 \\
  $\tau=10;\rho$=0; N=500&   9.94	&9.83	&2.31	&0.62 &  0.0057	&0.0051	&275.2240	&301.7934\\
      $\tau=10;\rho$=0.5; N=100   & 6.05&	5.94&	1.08	&0.63&  0.0262	&0.0126&	71.3873	&63.0084 \\
       $\tau=10;\rho$=0.5; N=500  &9.95&	9.82&	2.63&	0.67& 0.0056	&0.0050&	308.7230	&334.8119\\
             $\tau=10;\rho$=0.99; N=100 & 1.52&	1.26&	4.31&	2.26  & 0.0301	&0.0107&	1031.6216&	587.0757\\
       $\tau=10;\rho$=0.99; N=500 &   8.38	&7.56&	17.11	&3.84   &0.0046	&0.0041&	1315.9980&	1335.4355 \\
\hline\end{tabular}
}
    \caption{Summary of the true positives, false positives, root-mean-square prediction errors (RMSP), and computation time (seconds per repetition) for FSL and AFSL.  Here $N$ denotes the sample size, the number of points per curve is 50, $\tau$ indicates the strength of the within curve correlation for the errors, and $\rho$ is the strength of the between predictor correlation.}  \label{tab:111}
\end{table}
\clearpage

\begin{figure}
\centering
\includegraphics[width=0.78\textwidth, right]{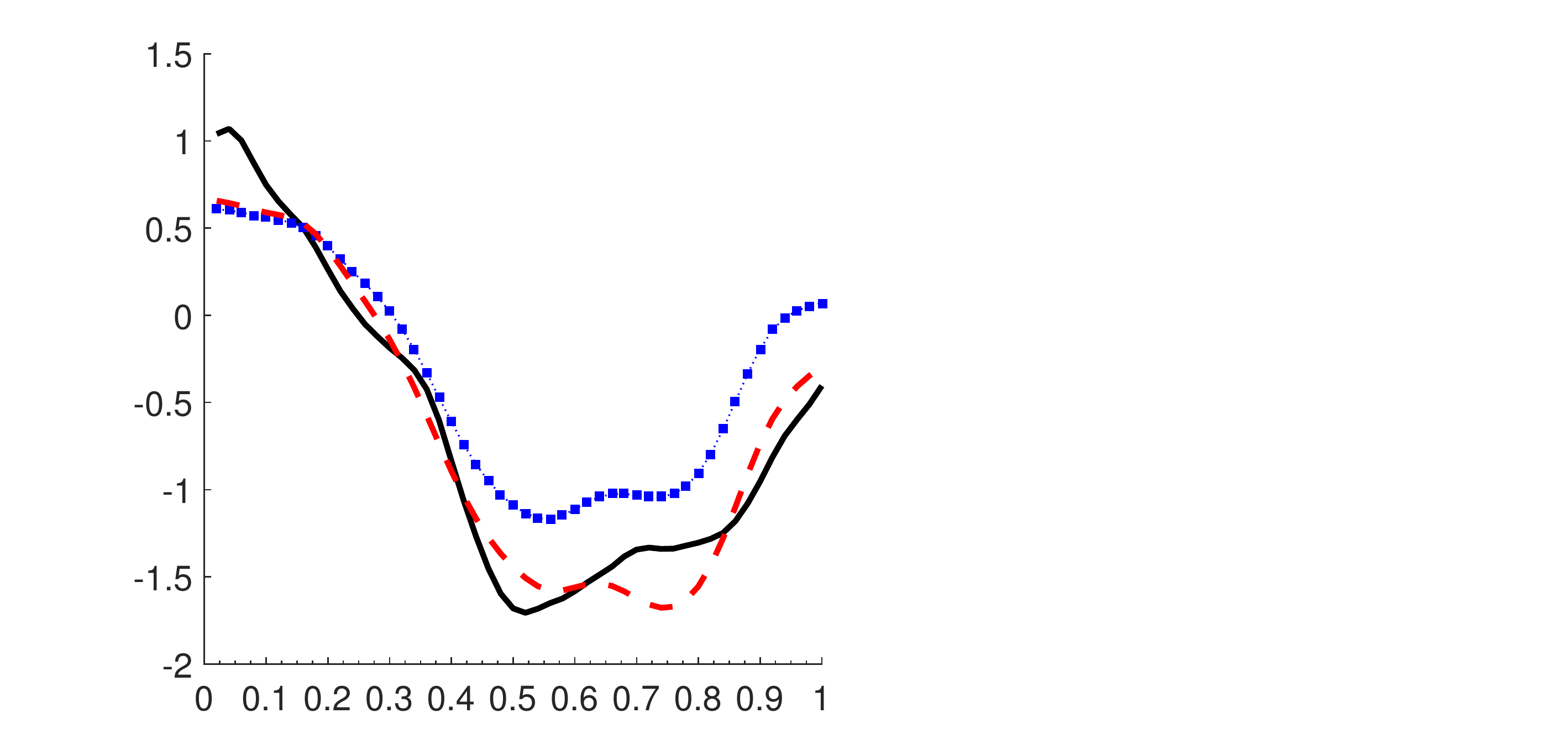}
\caption{ One example of an estimated beta versus a true one for $\tau_1$=0.25, $\tau_2=0.25$, $N = 500$ and $\rho =0$.The true $\beta$ is plotted using a solid black line, the FSL estimate is plotted using blue squares, and the AFSL estimate is plotted using a red dashed line.}
\label{tab:F1}
\end{figure}

\begin{table}[ht]
\centering
\resizebox{5.8in}{!}{%
\begin{tabular}{ccccccccl} \hline
& \multicolumn{2}{c}{True Positives} &  \multicolumn{2}{c}{False Positives}& \multicolumn{2}{c}{RMSP} &  \multicolumn{2}{c}{ Time (seconds)}\\ \hline
      Parameters & {FSL} & {AFSL}& FSL &  AFSL& {FSL} & {AFSL}& FSL &  AFSL\\ \hline
       $ \tau=0.01; \rho$=0; N=100 &  9.87	&9.82	&2.32&	0.42 & 0.0242&	0.0208	&118.892	&124.8117  \\
     $\tau=0.01;\rho$=0; N=500&9.99&	9.99&	6.98&	0.02&0.0092	&0.0091&	289.4841&	321.3352\\
      $\tau=0.01;\rho$=0.5; N=100   & 9.85	&9.66	&2.84&	0.55& 0.0246	&0.0218	&129.4541	&135.4117\\
       $\tau=0.01;\rho$=0.5; N=500 &9.99	&9.96&	6.34&	0.01&  0.0092	&0.0092	&306.2771	&335.4572\\
             $\tau=0.01;\rho$=0.99; N=100 &4.05&	3.57&	7.71	&3.39 & 0.0368	&0.0173&	856.7992	& 667.4319\\ 
       $\tau=0.01;\rho$=0.99; N=500 & 9.73	&9.34	&17.21&	0.33 &0.0071	&0.0068&	1153.8001&	1177.1012\\ \hline
$\tau=0.25; \rho$=0; N=100 & 8.35	&8.11&	1.37&	0.52 & 0.0348&	0.0235&	79.9695&	82.5338 \\ 
  $\tau=0.25;\rho$=0; N=500&9.99	&9.93	&2.45	&0.13&    0.0099	&0.0094	&317.6521	&344.4366\\ 
      $\tau=0.25;\rho$=0.5; N=100   &7.72	&7.55&	1.59	&0.72&0.0376	&0.0211	&  83.9127&	80.6943\\
       $\tau=0.25;\rho$=0.5; N=500 &9.99&	9.95&	2.41&	0.13&  0.0099&	0.0093	&284.5948	&307.7951\\
             $\tau=0.25;\rho$=0.99; N=100 &1.38	&1.21& 3.19&	1.84&0.0551	&0.0149	&921.1164	&433.9953\\ 
       $\tau=0.25;\rho$=0.99; N=500 &8.83	&8.25	&15.83	&2.86&0.0082&	0.0073	&1114.5110&	1131.1231\\ \hline
$\tau=1; \rho$=0; N=100 & 6.68&	6.57	&1.13	&0.62&  0.0424&	0.0195	&61.2572	&56.2244 \\  
  $\tau=1;\rho$=0; N=500& 9.99&	9.91	&2.42	&0.45&    0.0103	&0.0094	&285.0516	&310.0837\\
      $\tau=1;\rho$=0.5; N=100   & 6.83	&6.68&	1.55	&0.87&   0.0417	&0.0199	&79.4073&	72.5188 \\
       $\tau=1;\rho$=0.5; N=500 &9.98	&9.91	&2.25	&0.42&   0.0103	&0.0095&	274.4325&	298.9810\\
             $\tau=1;\rho$=0.99; N=100 & 1.08	&0.92&	2.49	&1.52& 0.0581&	0.0145	&926.1381&	392.3812\\ 
       $\tau=1;\rho$=0.99; N=500 & 8.80&	8.06	&16.25	&3.39 &0.0083	&0.0074&	 1129.6421&	1145.9041\\ \hline
$\tau=10; \rho$=0; N=100 &   6.15&	6.04&	0.72	&0.47  &  0.0455&	0.0196	&57.8118	&50.4278 \\
  $\tau=10;\rho$=0; N=500&  9.94	&9.91	&2.57&	0.62&  0.0103	&0.0094	&248.1699	&271.3013\\ 
      $\tau=10;\rho$=0.5; N=100   & 5.61&	5.54	&0.97	&0.57& 0.0477	&0.0190&	64.7757	&53.0382 \\
       $\tau=10;\rho$=0.5; N=500  &9.97&	9.89&	2.69	&0.58&0.0104	&0.0096	&257.8873	&280.6499\\
             $\tau=10;\rho$=0.99; N=100 &1.12&	1.02	&2.54	&1.41 &  0.0578&	0.0141	& 988.6308 &403.8979\\
       $\tau=10;\rho$=0.99; N=500 & 9.01	&8.18&	15.92&	3.42&0.0083&	0.0074&	1114.6311&	1133.2537\\  
\hline\end{tabular}
}
    \caption{Summary of the true positives, false positives, root-mean-square prediction errors (RMSP), and computation time (seconds per repetition) for FSL and AFSL.  Here $N$ denotes the sample size, the number of points per curve is 16, $\tau$ indicates the strength of the within curve correlation for the errors, and $\rho$ is the strength of the between predictor correlation.}  \label{tab:222}
\end{table}
\clearpage
 \textcolor{black}{Turning first to variable selection consistency in Table \ref{tab:111}, we see that AFSL cannot beat FSL in terms of selecting true predictors.  This is because AFSL uses the FSL selection as a starting point.  The real gain is in terms of the false positives.  For low and moderate correlation, we see that the AFSL has about $50\%$ lower false positive rate for $N=100$, and this increases dramatically for $N =500$ resulting in an over $90\% $ decrease in the false positive rate.  This is accomplished while maintaining nearly exactly the same true positive selection rate.  Both methods, as we would expect, have a noticeably harder time when the correlation between predictors is high, but AFSL still substantially decreases the number of false positives.
Turning to the prediction error (RMSP) the differences are not quite as dramatic, but AFSL is still out performing FSL.  For smaller samples, the gain is around 10-20\%, but this gain reduces for larger $N$.  Lastly, adding AFSL on top of FSL increases the computation time by less than 10\%. We therefore conclude that there is very little reason not to run AFSL after FSL. To help illustrate the resulting estimates, in Figure \ref{tab:F1} we plot one example of an estimated beta versus a true one for $\tau_1$=0.25, $\tau_2=0.25$, $N = 500$ and $\rho =0$.  There we can better visualize how AFSL compares to FSL.  In general, FSL will pull the estimates to zero more than AFSL.  This results in too much bias for the nonzero betas, and thus worse statistical performance.}

 \textcolor{black}{As our last set of simulations, we rerun all of the scenarios in Table \ref{tab:111}, but reduce the number of time points sampled per curve to 16 to better emulate the CAMP data.  The results are reported in Table \ref{tab:222}.  As we can see, the performance of the two methods is nearly the same as before.  This helps justify  applying our methodology to the CAMP data, which we will discuss more in the next section.}

\subsection{Data Example}
\begin{figure}
\centering
  \centering
  \includegraphics[width=0.6\linewidth]{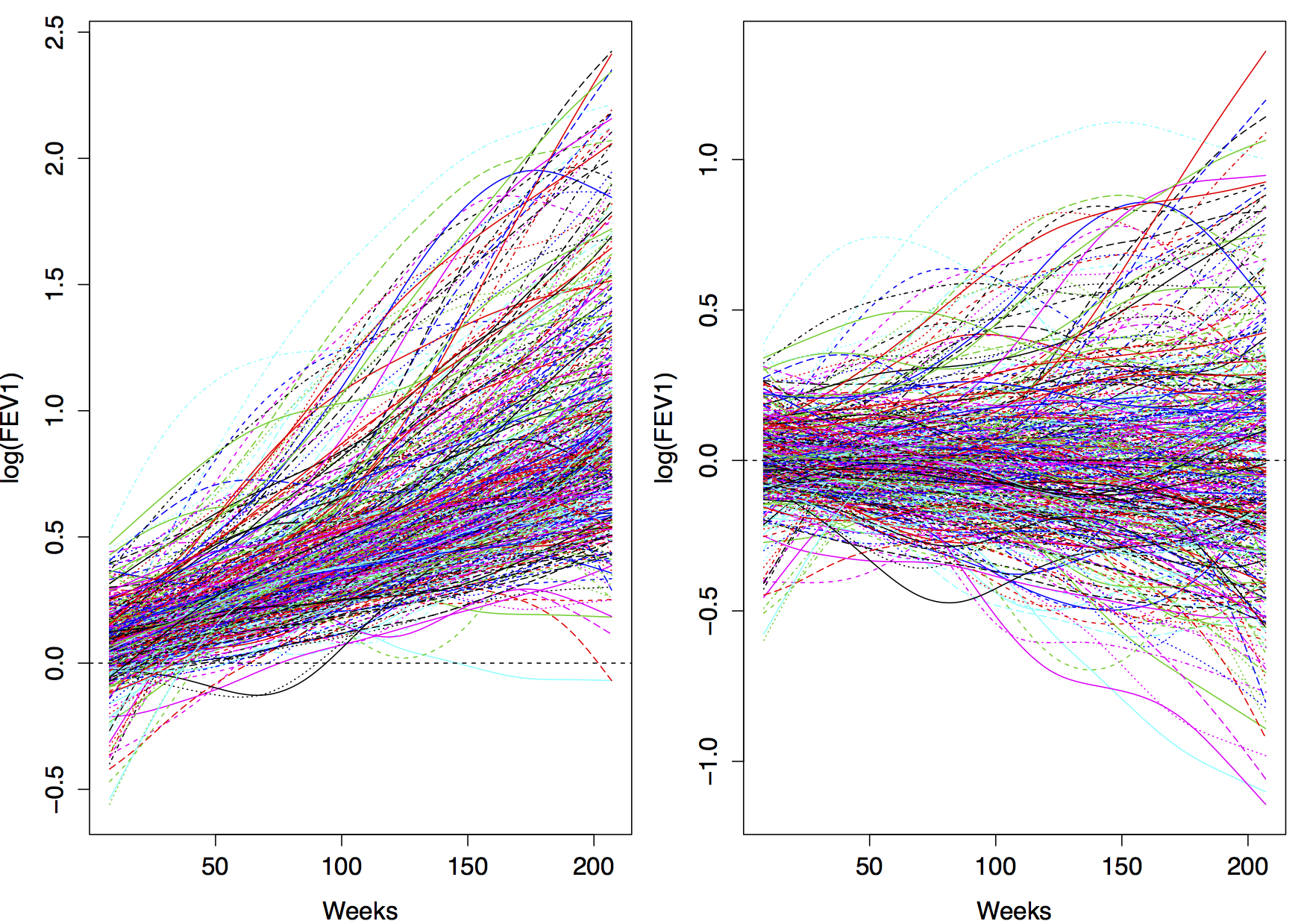}
\caption{\textcolor{black}{The left panel is a plot of the log(FEV1) curves for the 540 subjects from CAMP. Each trajectory  is estimated using penalized bsplines.  The right panel is the same, but the curves have been standardized by age, gender, and treatment.   }} \label{tab:55}
\end{figure}

In this section we apply FSL and AFSL to data from the Childhood Asthma Management Program, CAMP.  CAMP was a multi--center, longitudinal clinical trial designed to better understand the long term impact of two common daily asthma medications, Budesonide and Nedocromil, on children \cite{camp:1999}.  Phenotypic and genotypic data are freely available at dbGaP, study accession phs000166.v2.p1.  The data we consider here consists of 540 caucasian subjects, ages 5--12, with asthma, who were randomly assigned a particular treatment (one of the two drugs or placebo).  Each subject made 16 clinical visits spread over a four year period. The outcome is (log) \textit{forced expiratory volume in one second} or FEV1, which is a commonly used proxy for lung strength.  For the domain, we use ``study time" which represents how long the subject has been on a particular treatment.  In  Figure \ref{tab:55} \textcolor{black}{we plot log(FEV1) curves for the 540 subjects. The left panel shows curves expanded using bsplines, with the same approach as in Section \ref{s:sims}.  The right panel plots the curves after regressing out gender, age, and treatment, using standard function-on-scalar regression techniques \cite{re14,re15}.  We use the standardized curves as response functions when applying our methods.}

The data are observed on a grid, which is not quite evenly spaced as there are more observations early on in the study.  In particular, the first three observations are 1-2 weeks apart, the next two are 2 months apart, and then remaining 11 visits are 4-5 months apart.  As with the simulations, the data are preprocessed using bsplines and FPCA.  We mention that with 16 observed time points, readers will rightly be concerned about the effects of such preprocessing on data which is not densely observed.  However, since the data is observed on a common grid and 16 is well beyond the typical rule of thumb ($540^{1/4} \approx 5) $ we believe that this approach is justified.  Furthermore, as we saw in the simulations Section \ref{s:sims}, our methods still work well with 16 time points. 

\begin{figure}
\centering
\begin{subfigure}{.248\textwidth}
  \centering
    \caption*{rs7751381}
  \includegraphics[width=2.1\linewidth]{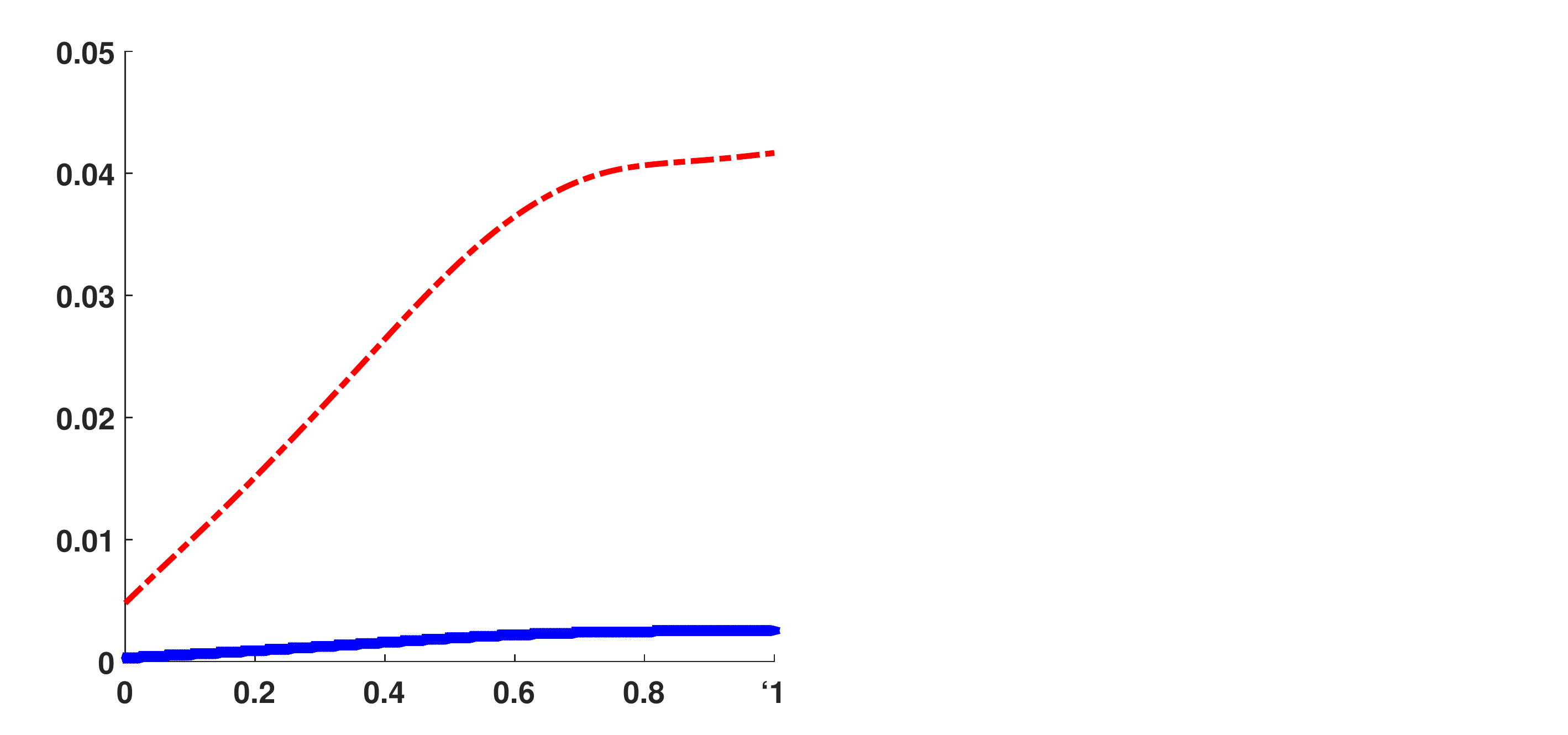}
  \end{subfigure}%
\begin{subfigure}{.248\textwidth}
  \centering
    \caption*{rs1540897}
  \includegraphics[width=2\linewidth]{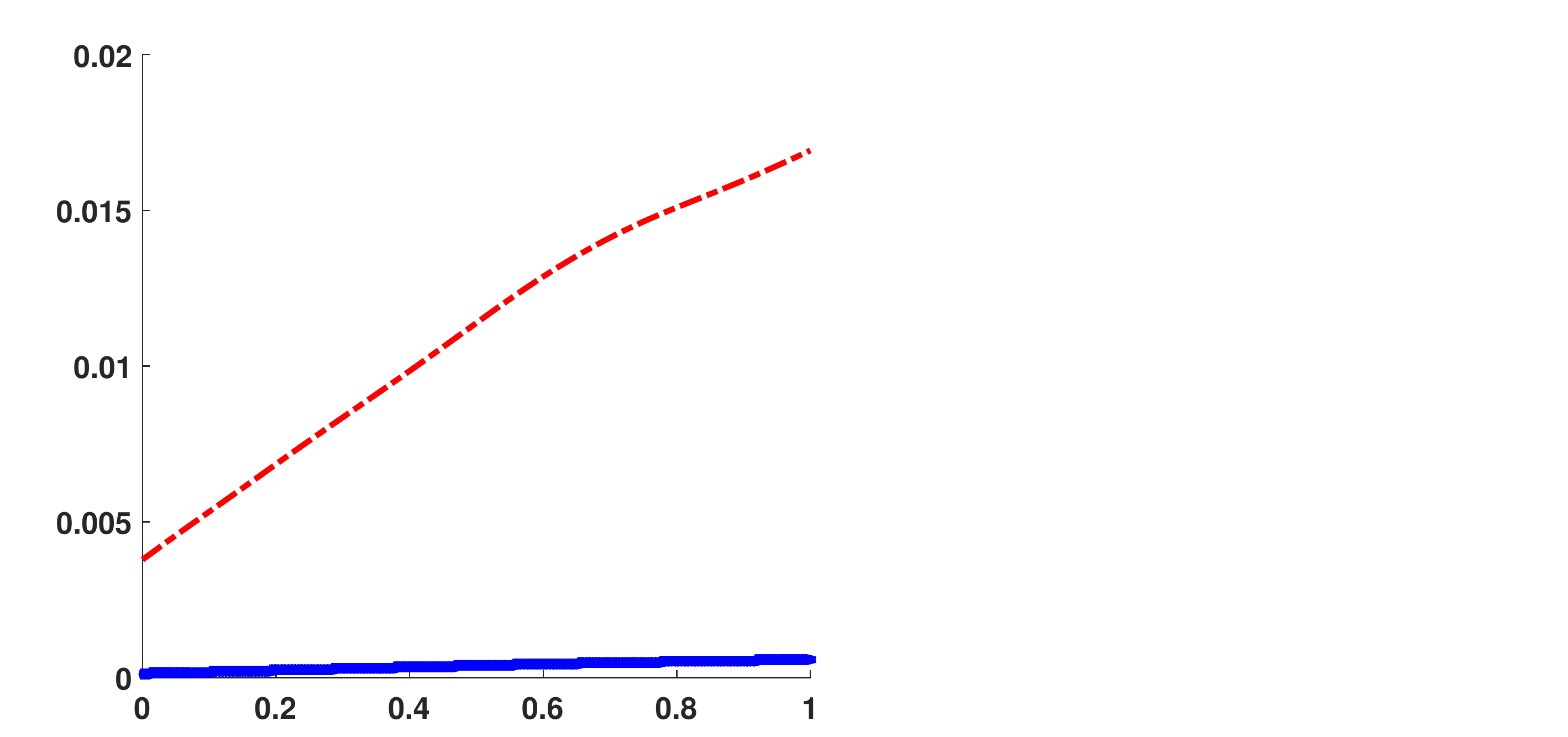}

\end{subfigure}
\begin{subfigure}{.248\textwidth}
  \centering
    \caption*{rs17372029}
  \includegraphics[width=2.1\linewidth]{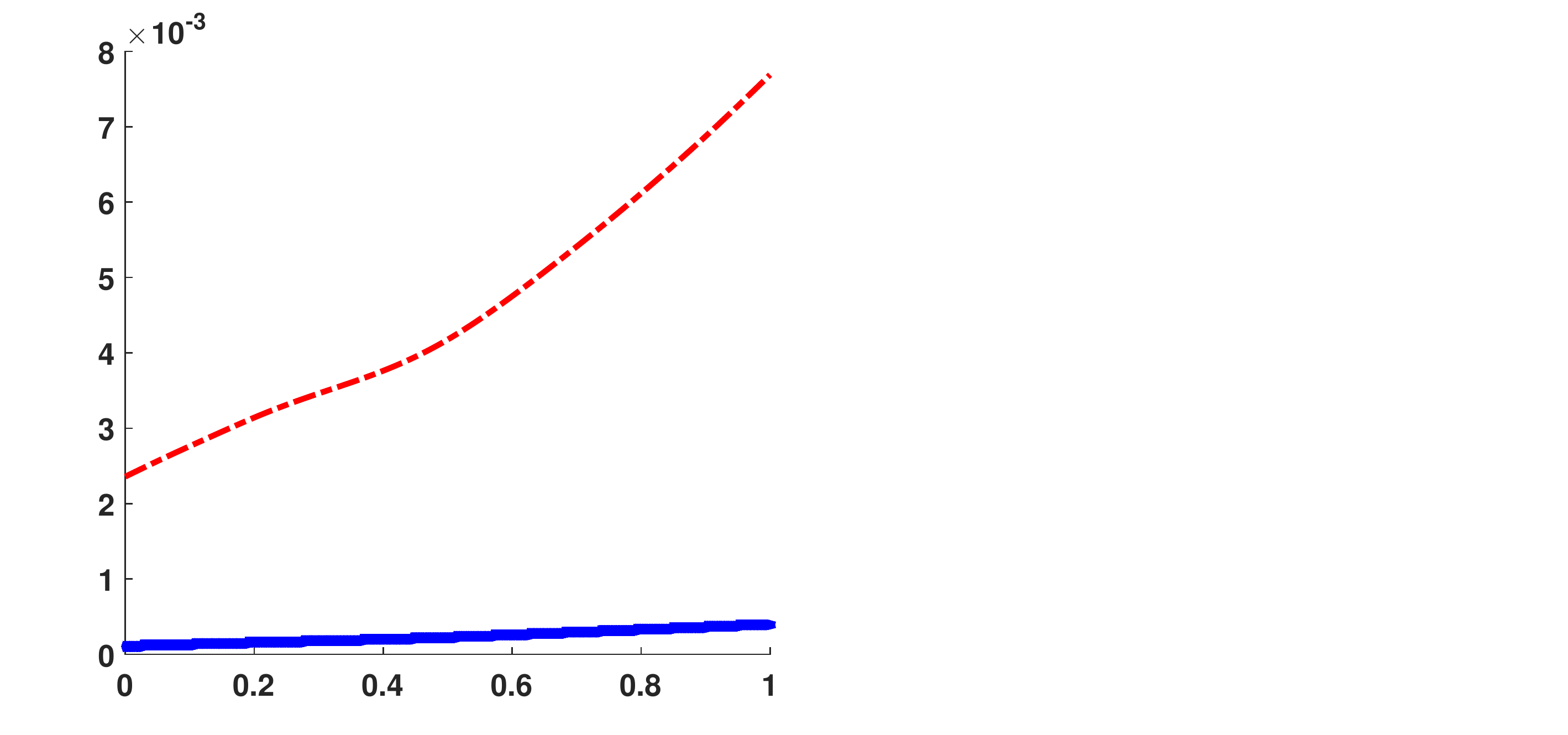}

\end{subfigure}%
\begin{subfigure}{.248\textwidth}
  \centering
  \caption*{rs1875650}
  \includegraphics[width=2\linewidth]{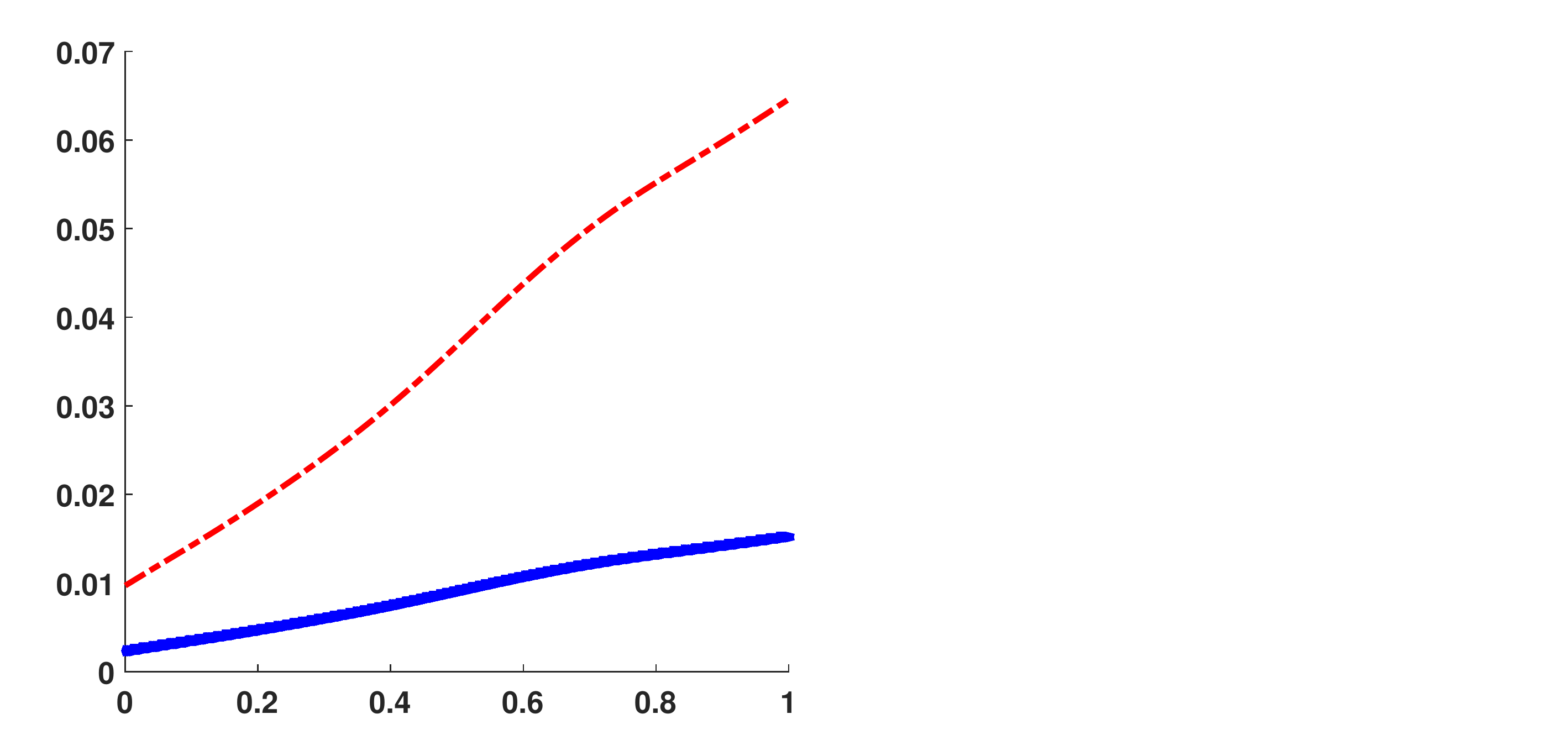}

\end{subfigure}

\begin{subfigure}{.248\textwidth}
  \centering
    \caption*{rs2019435}
  \includegraphics[width=2\linewidth]{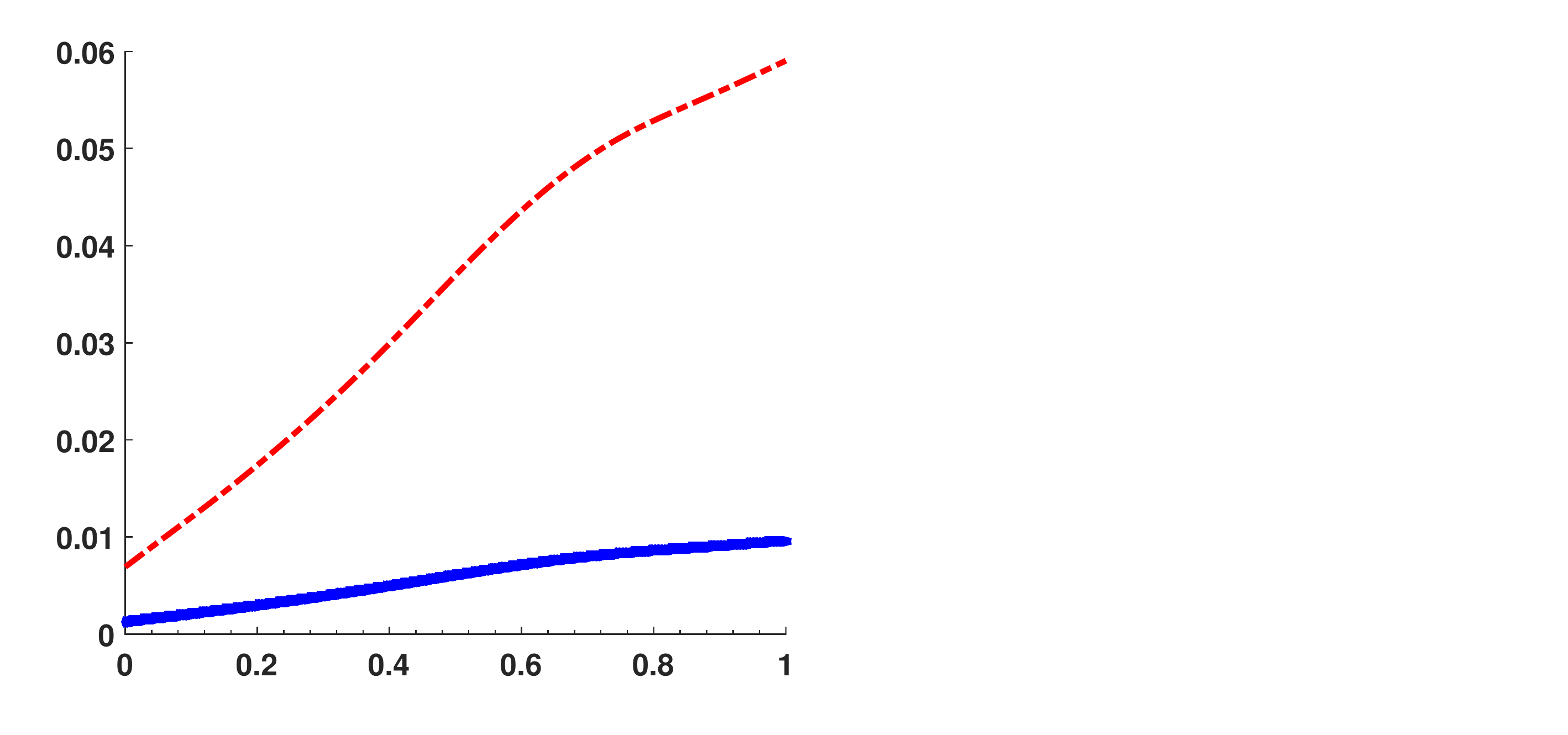}

\end{subfigure}%
\begin{subfigure}{.248\textwidth}
  \centering
      \caption*{rs4734250}
  \includegraphics[width=2\linewidth]{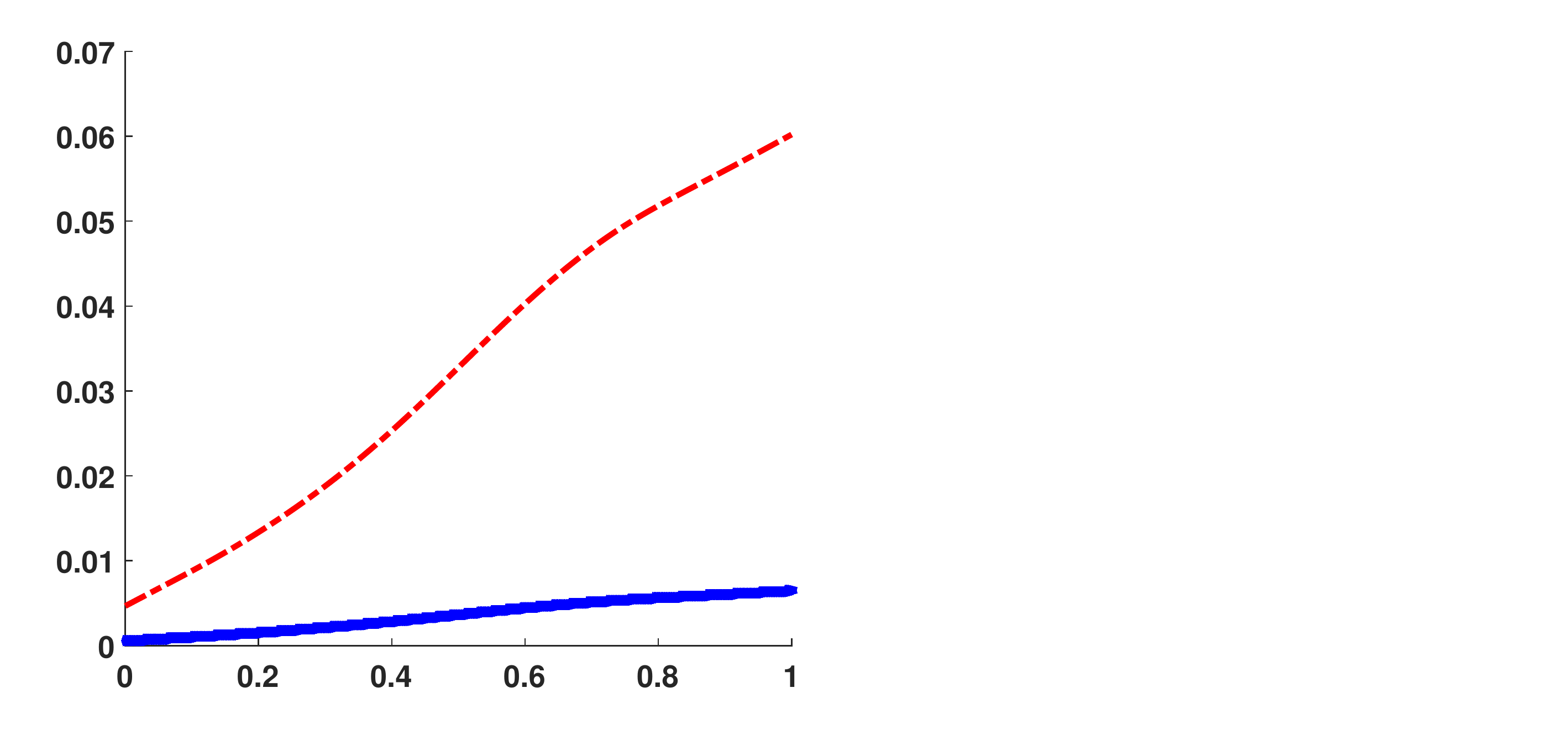}

\end{subfigure}
\begin{subfigure}{.248\textwidth}
  \centering
    \caption*{rs4752250}
  \includegraphics[width=2\linewidth]{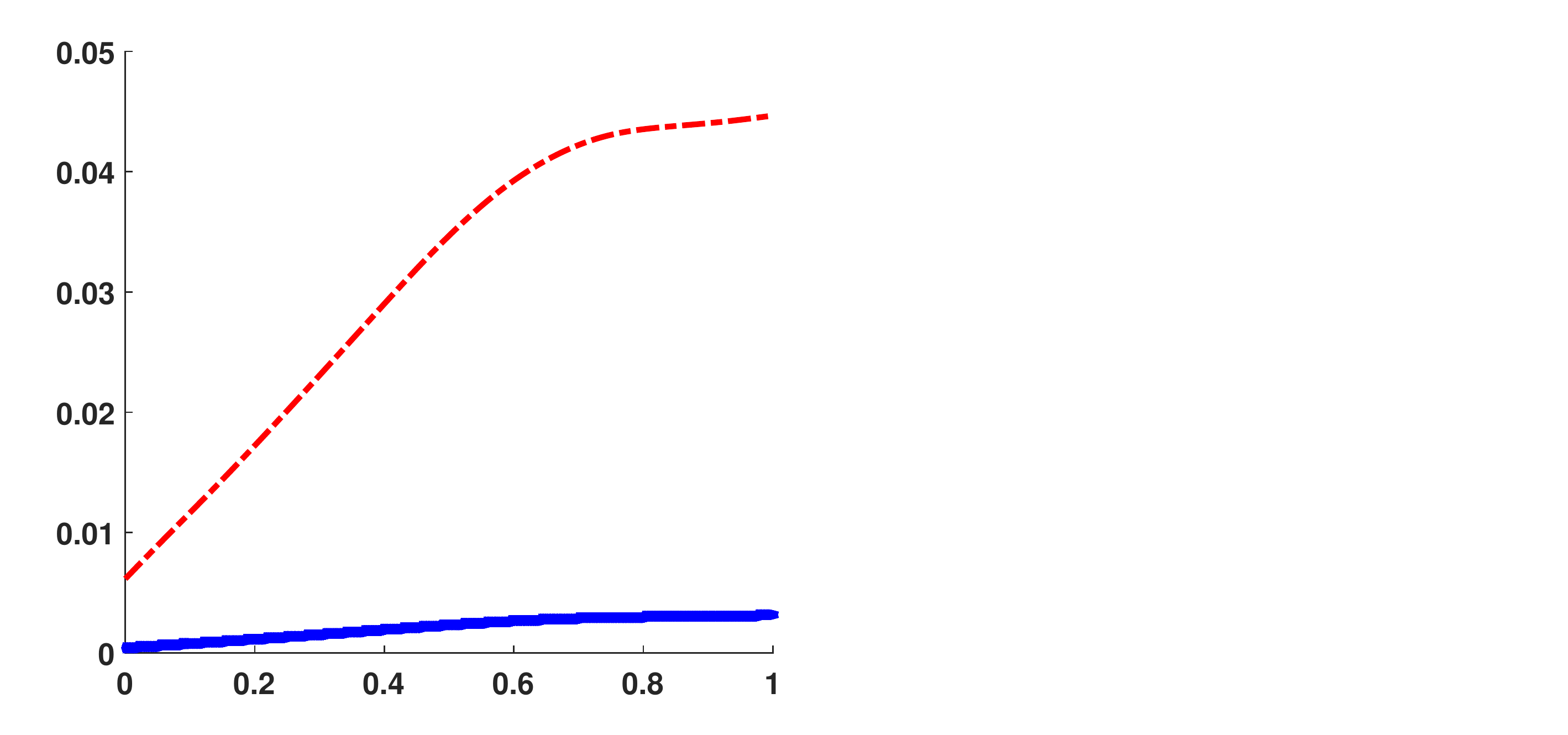}
\end{subfigure}%
\begin{subfigure}{.248\textwidth}
  \centering

     \caption*{rs2206980}
  \includegraphics[width=2.1\linewidth]{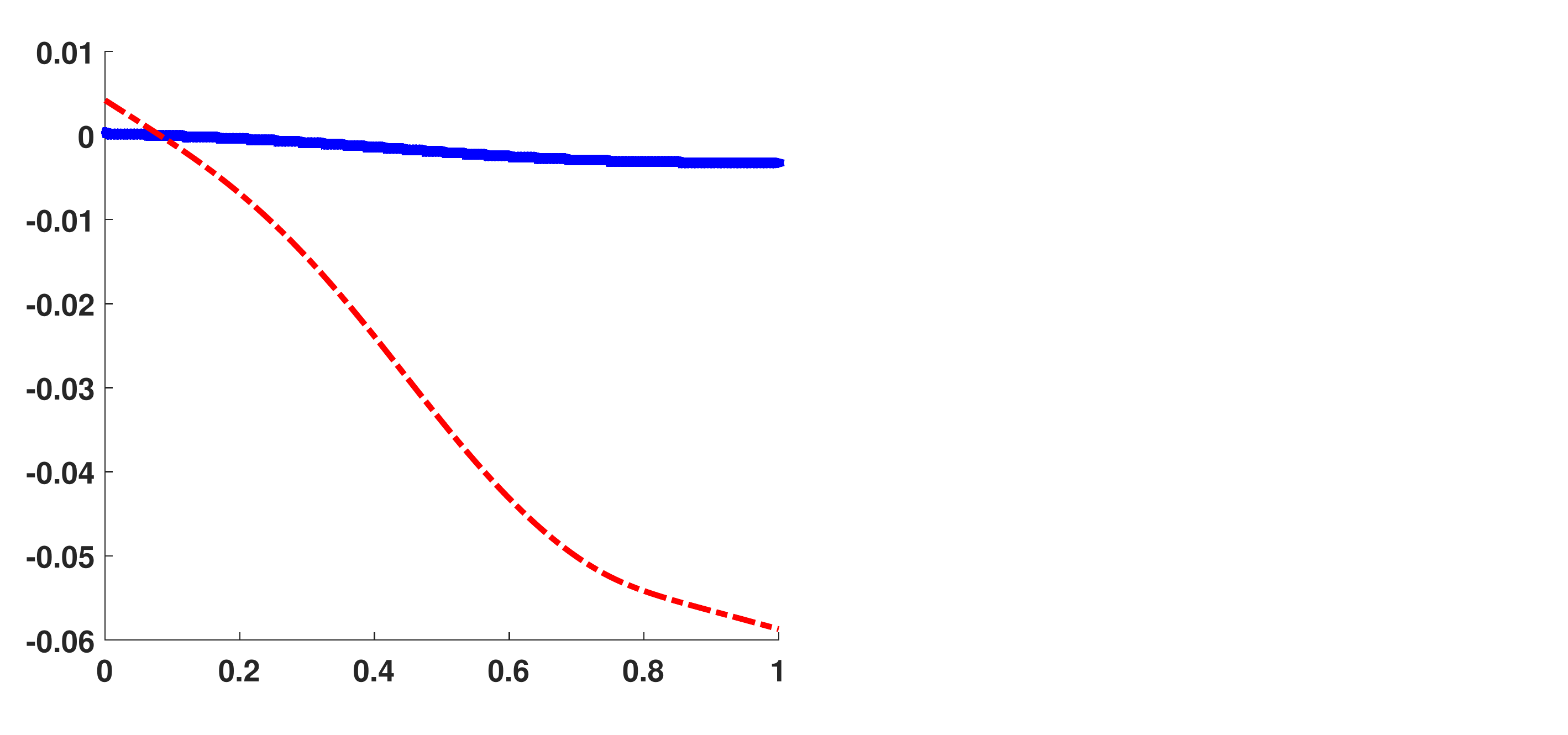}
\end{subfigure}

\begin{subfigure}{.248\textwidth}
  \centering
   \caption*{rs1368183}
  \includegraphics[width=2\linewidth]{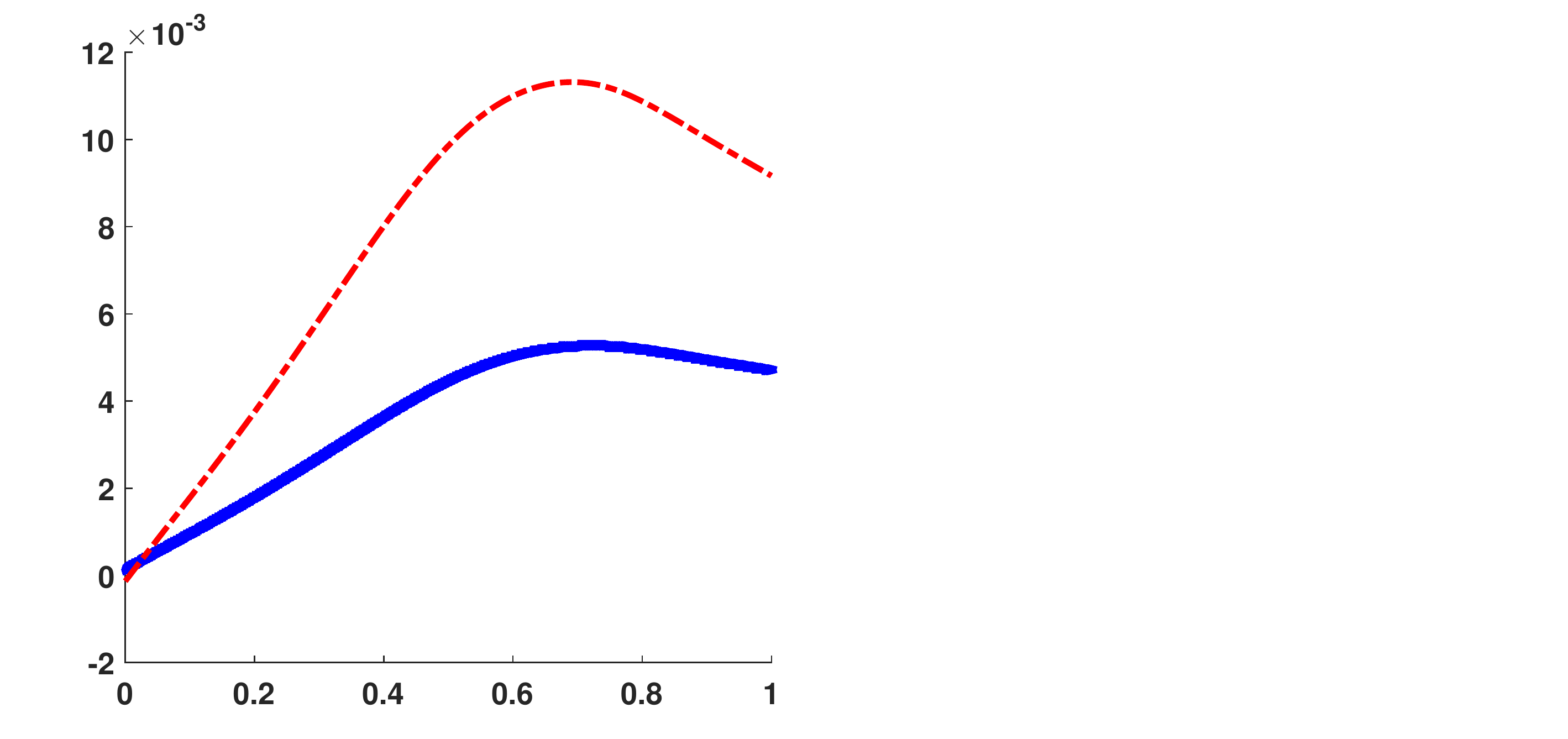}

  \label{fig:sub1}
\end{subfigure}%
\begin{subfigure}{.248\textwidth}
  \centering
   
      \caption*{rs953044}
  \includegraphics[width=2.1\linewidth]{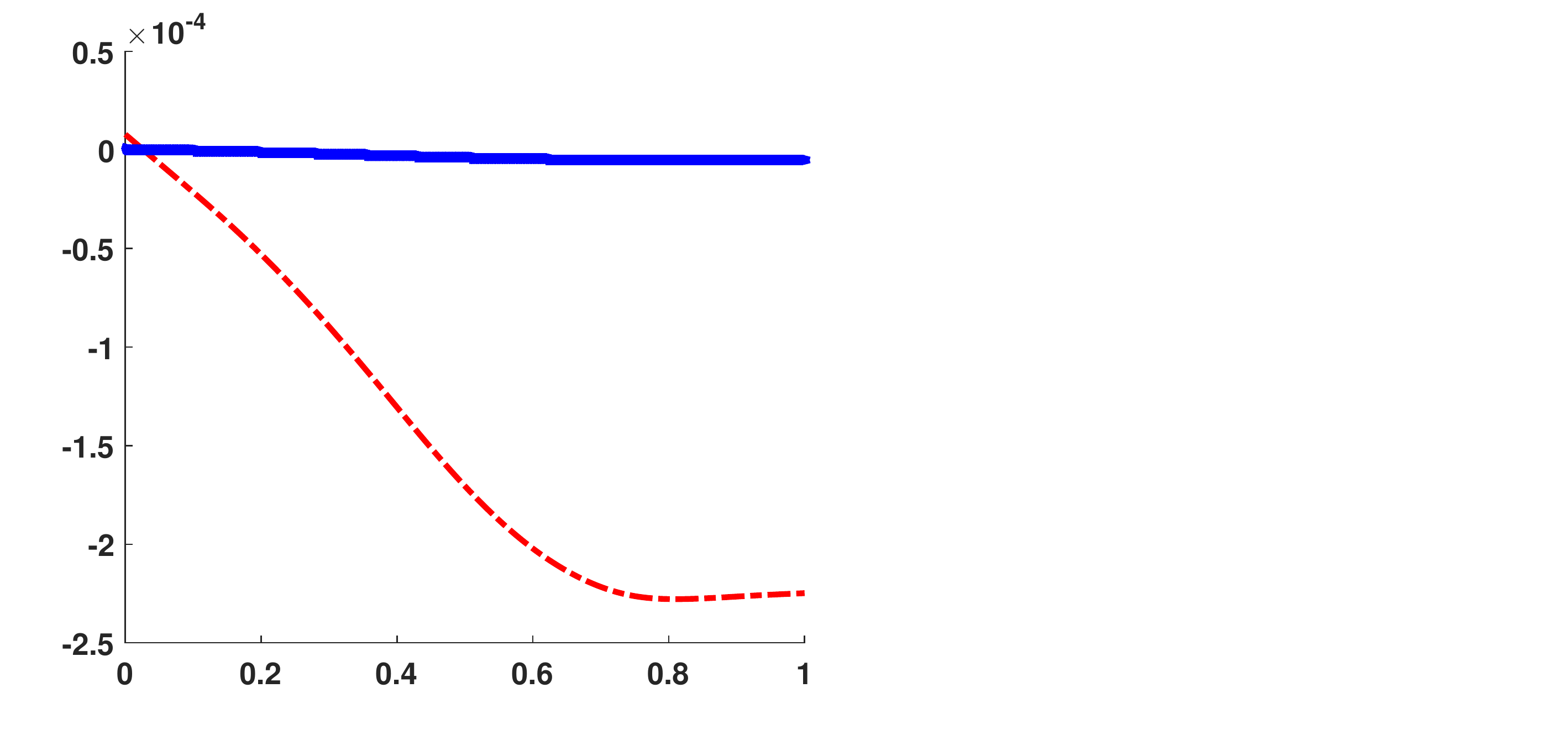}

  \label{fig:sub22222}
\end{subfigure}

\caption{Estimated $\beta$ coefficients for top ten SNPs from CAMP. The FSL estimate is plotted using blue squares and the AFSL is plotted using red dashes.}\label{fig:test}
\end{figure}

We applied FSL and AFSL to a subset of 10,000 of single nucleotide polymorphisms, SNPs, which were prescreened using methods from  \citetext{Ch:Li:Re:2016}.  
As in \citetext{M1}, we applied FSL and AFSL using both BIC and extended BIC (with a parameter of $0.2$).  It is well known that the extended BIC helps control false positives better when the number of predictors is large \cite{chen:2008,foygel:drton:2010}.  In Table \ref{t:snps} 11 out of 10,000 SNPs were selected by FSL, and 10 out of 11 SNPs selected by AFSL, regardless of whether BIC or EBIC was used. 
In Figure \ref{fig:test}, we see that rs1875650, rs1368183, rs7751381, rs17372029, rs1540897, rs4734250, rs4752250 and rs2019435 have a positive impact on lung development in growing children. We can also conclude that
rs953044 and rs2206980 have a negative impact on lung development in growing children.  Comparing FSL with AFSL, we see that FSL tends to over shrink the estimates and a great deal of the curvature of the estimates is lost.  In contrast, AFSL shrinks less and the changes over time are much clearer.  \\

\begin{table}[ht]
\begin{centering}
\resizebox{4in}{!}{%
\begin{tabular}{llllll}
\hline
 \multicolumn{2}{c}{SNP} & \multicolumn{2}{c}{FSL}  & \multicolumn{2}{c}{AFSL}    \\ \hline
 Chromosome & Name &BIC  & EBIC    &BIC  & EBIC \\ \hline
1 & rs1875650& $\times$ &$\times$&$\times$& $\times$ \\\hline
 2 & rs953044 & $\times$ &$\times$  & $\times$& $\times$         \\ \hline
 5 & rs1368183& $\times$ &$\times$  & $\times$& $\times$   \\ \hline
 6 & rs7751381& $\times$ &$\times$  & $\times$& $\times$ \\ \hline
 6 & rs2206980 & $\times$ &$\times$  & $\times$& $\times$  \\ \hline
 7 & rs17372029& $\times$ &$\times$  & $\times$& $\times$    \\ \hline
 8 & rs1540897& $\times$ &$\times$  & $\times$& $\times$ \\ \hline
 8 & rs4734250  & $\times$ &$\times$  & $\times$& $\times$        \\ \hline
10 & rs4752250& $\times$ &$\times$  & $\times$& $\times$  \\ \hline
 15 & rs2019435& $\times$ &$\times$  & $\times$& $\times$ \\ \hline
 20 & rs2041420& $\times$ &$\times$  && \\ \hline
\end{tabular}
}
\caption{Top SNPs selected by FSL and by AFSL. A $\times$ indicated the SNP is selected when using the corresponding variable selection criteria.}\label{t:snps}
\end{centering}
\end{table}

\textcolor{black}{
As a final comparison we compare the prediction error between FSL and AFSL on the CAMP data.  To do this, we split the data set into a training set (80$\%$) and test set (20$\%$). 
We build our model on the training set while comparing the predictive performance on the test set.  We measure performance in terms of the root mean squared error, RMSE:\begin{equation}
    \textrm{RMSE} = \sqrt{\frac{1}{N} \sum_{n=1}^{N} \|Y_n - \hat{Y}_n\|^2}.
\end{equation}
Here $\hat Y_n$ is the predicted functional outcome and $Y_n$ is the observed for the $n$th subject from the test dataset. 
We repeat this split 10 times and average the results, which are given in Table \ref{tab:11}.  
The difference in terms of BIC vs EBIC seems to be negligible, while AFSL consistently produces an RMSE which is about $15\%$ lower than FSL.  This further supports the idea that rs20141420 might be a false positive, as this was the only one dropped by AFSL, as indicated in Table \ref{t:snps}}.

\begin{table}[ht]
\centering
\resizebox{3in}{!}{%
\label{my-label}
\begin{tabular}{lllll}
\hline
& \multicolumn{2}{c}{BIC}  & \multicolumn{2}{c}{EBIC}    \\ \hline

Test   & FSL & AFSL    & FSL & AFSL  \\ \hline
1  &        0.1472  &      0.1310 &       0.1386&	0.1194\\ \hline
2  &   0.1533  & 0.1275    &   0.1536&	0.1272\\ \hline
3  &    0.1228
   &  0.1267    &   0.1523&	0.1062 \\ \hline
4  &    0.1429  &    0.1098 &   0.1444&	0.1315\\ \hline
5  &  0.1257  &       0.1299&  0.1372&	0.1248   \\ \hline
6  &   
   0.1418  &  0.1178  & 0.1431&	0.1167  \\ \hline
7  &       0.1576  &  
    0.1219   &      0.1372&	0.1248  \\ \hline
8  &    0.1543  &      0.1094 &    0.1536&	0.1272   \\ \hline
9  &      0.1429 &  0.1098   &      0.1264&	0.1338
     \\ \hline
10 &     0.1320
   &    0.1502 &     0.1536	&0.1248 \\ \hline \hline 
Average &     0.1421	&0.1234  &0.1440	&0.1236\\ \hline 
\end{tabular}
}
\caption{ RMSE for the CAMP data calculated by FSL and AFSL when using the corresponding variable selection criteria.}   \label{tab:11}
\end{table}

\section{Conclusion and Future Work}
In this paper we have presented an adaptive version of the function-on-scalar lasso introduced in \citetext{M1}.  Using techniques from convex analysis over Hilbert spaces, we were showed that AFSL estimates achieve the functional oracle property even when the number of predictors is exponential relative to the sample size.  

There are still a great deal of open problems in this area as relatively little has been done for high dimensional functional data.  We used AFSL to control the bias inherent in FSL, but other approaches should work as well.  For example, SCAD would also be an excellent candidate for reducing the bias of the estimates, but serious theoretical developments would be required.  In another direction, alternating the penalty to more carefully control the level of smoothing could prove beneficial as well.  
Additionally, these methods would greatly benefit from customized computational tools.  In the present work and in \citetext{M1}, an ADMM for scalar group lasso was used to find the FSL and AFSL estimates.  While this approach works well, it can be improved by producing custom ADMM or coordinate descent methods which are specifically designed for functional data.

Lastly, the tuning parameter was selected using BIC and EBIC, and in practice they work well, but a stronger theoretical justification for their use is still needed in the functional context.  It would also be useful to explore other tuning parameter methods and see how they compare.

\clearpage
\appendix
\section{Proofs}
In this section we collect the proofs for Theorems \ref{t:subgrad} and \ref{t:oracle}.  The former is fairly straightforward as the structure is nearly the same as for the scalar or multivariate setting.  The latter is much more involved and will be broken into several pieces.

\subsection{Proof of Theorem \ref{t:subgrad}}
For a differentiable convex function, the subdifferential is equal to the (Fr\'echet) derivative.  So we have that the subdifferential of the first term of $L_\lambda^A$ with respect to $\beta_i$  is given by
\[
- \frac{1}{N}\sum_{n=1}^N X_{ni}(Y_n - X_n^\top \beta).
\]
Turning to the second term, consider the functional $g(x) = \| x\|$.  For $x \neq 0$, we need to verify that subgradient is $\|x \|^{-1} x$, i.e., that for all $y$
\begin{align*}
g(y) - g(x) = \|y\| - \|x\| \geq \|x \|^{-1} \langle x, y - x \rangle.
\end{align*}
This is equivalent to showing that
\[
 \|y\|\|x\| -  \|x\|^2  \geq  \langle x,y\rangle- \|x\|^2,
\]
but this follows from the Cauchy-Schwarz inequality.  For $x \neq 0$ we need to show that for any $\|h \| =1$ and all $y$ we have
\[
g(y) - g(x) = \|y\|  \geq  \langle h, y \rangle,
\]
which again just follows from Cauchy-Schwarz.  Thus the subdifferential of $g(x)$ is verified.  Combining this fact with the additivity and linearity of subdifferentials, we have that Theorem \ref{t:subgrad} holds.

\subsection{Theorem \ref{t:oracle}: Preliminary}
Here we state several key preliminary results which are useful for the main proof.  We begin with a reminder of the restricted eigenvalue condition.
\begin{definition} (Restricted eigenvalue condition). A matrix $M \in \mathbb{R}^{N\times I}$ satisfies the RE$(I_{0},\alpha)$ condition if, for all subsets S $\subset \{1,\cdots, I\}$ with $|S|\leq I_{0}$, 
\begin{equation*}
\|Mw\|_{2}^{2}\geq \alpha N\|w\|_{2}^{2} ~~ \text{for all} ~~w \in \mathbb{R}^{I}~~\text{with}~~\|w_{s^{c}}\|_{l_{1}}\leq 3\|w_{s}|_{l1}
\end{equation*}
\end{definition}

In \citetext{M1} we showed that the scalar restricted eigenvalue condition implied a functional one which was necessary for the convergence rates in FSL.  We now restate a result which follows from \citetext{van:2009}, namely, that the irreprsentable condition combined with a sparsity assumption implies that the restricted eigenvalue condition holds.
\begin{lemma}
If a matrix $\bM$ satisfies  Assumption \ref{a:tech} then it also satisfies RE$(I_{0},\alpha)$.
\end{lemma}
Given the above result we can apply Theorem 3 from \citetext{M1} to obtain the following.	
\begin{lemma} \label{l:fsl}
If Assumptions \ref{a:model} and \ref{a:tech} hold then the FSL estimate $\tilde \beta$ satisfies 
$$\sup_{i\in \mathcal{S}} \|\beta_{i}-\tilde \beta_{i}\|=O_{p}(r_{N}^{1/2})
\quad \text{where}
\quad r_N = \frac{I_0 \log(I)}{N }.
$$
\end{lemma}

The following functional concentration inequality from \citetext{M1} will also be useful.
\begin{lemma}\label{l:barb}
 Let $X$ be an $\mathcal{H}$ valued Gaussian process with mean zero and covariance operator
C. 
Then we have the bound \\
		\begin{equation*}
	 P\{\|X\|^{2}\geq \|C\|_{1}+2\|C\|_{2}\sqrt{t}+2\|C\|_{\infty}t \}\leq \exp(-t).
		\end{equation*}
		where $\|C\|_{1}$ is the sum of the eigenvalues of $C$, $\|C\|_{2}^2$ is the sum of the squared eigenvalues, and $\|C\|_{\infty}$ is the largest eigenvalue.
\end{lemma}

Our last lemma is mainly for technical convenience as it shows that the operator norm of a matrix is the same regardless of the Hilbert space it is acting on.  This is especially useful as it implies that the functional and scalar irrepresentable conditions are the same.
\begin{lemma}
Let $ \bA \in \mbR^{q \times p}$ be viewed as a linear operator from $ \mcH ^{p}\to \mcH^{q}$.  Then the operator norm of $\bA$ is the same for any Hilbert space $\mcH$, and in particular when $\mcH = \mbR$.
\begin{proof} 
Recall that  when $\bA$ is viewed as $\bA: \mbR^{p}\to \mbR^{q}$. Then, 
\begin{align*}
\|\bA\|_{op}&=\sup_{|x|=1}|\bA x |= \sigma_1,
\end{align*}
where $\sigma_1$ is the largest singular value.	
When $\bA$ is viewed as on operator on another Hilbert space $\bA: \mcH^{p}\to \mcH^{q}$ we define
\begin{align*}
\|\bA\|_{F; op}&=\sup_{x\in\mcH^{p }: \|x\|=1}\|\bA x\|.
\end{align*}
Let $\{e_j\}$ be an orthonormal basis of $\mcH$, then $x_i = \sum_j x_{ij} e_j$, where $x_{ij} \in \mbR$.  Define $x^j = \{x_{ij}: i =1 \dots p\}$.  Then by Parseval's identity
\[
\| \bA x\|^2 = \sum_{j} | \bA x^j|^2 \leq \sigma_1^2 \sum_{j} |x^j|^2  = \sigma_1^2 \| x\|^2.
\]
Which means that $\| \bA\|_{F,op} \leq \| \bA\|_{op}$.  To see that they are equal, take $x$ such that $x^{1}$ is equal to the largest right singular vector of $\bA$, and all other $x^j$ are zero.  Then we have $\| \bA x\|^2 =  | \bA x^1|^2 = \sigma_1^2 = \| \bA\|_{op}^2 $, thus the two operator norms are equal.
\end{proof}
\end{lemma}

\subsection{Theorem \ref{t:oracle}: Proof}
We now have the necessary tools to establish the proof of Theorem \ref{t:oracle}.  We will split this into multiple pieces.  We begin with the longer argument, namely, that AFSL is variable selection consistent.  We will then finish with the shorter argument, namely, that AFSL is also asymptotically equivalent to the oracle estimate. \\

\noindent {\bf Proof of Theorem \ref{t:oracle}.2:} \\
We will mimic the structure of \citetext{h08a}.  Partition $\hat \beta  = \{\hat \beta_1, \hat\beta_2\}$, where $\hat \beta_1$ is the estimate of the true predictors, and $\hat \beta_2$ of the null predictors.  Define
\[
\tilde s_{1} = \{ \tilde w_{i} \| \hat \beta_{i} \|^{-1} \hat \beta_{i}, i \in \mcS\}.
\]
Suppose that $\hat \mcS = \mcS$, that is, the support of $\hat \beta$ is $\mcS$.  Then from Corollary \ref{c:sub} we have that
\[
\bX_1^\top (Y - \bX_1 \hat \beta_1) = \lambda \tilde s_1.
\]
	Solving for $\hat \beta_1$ we get that
	\[
	\hat \beta_1 = (\bX_1^\top \bX_1)^{-1}(\bX_1 Y - \lambda \tilde s_1)
	\]
	or equivalently
	\[
	\hat \beta_1 = (\bX_1^\top \bX_1)^{-1} (\bX_1^\top Y - \lambda \tilde s_{1})= \beta_1 +  (\bX_1^\top \bX_1)^{-1}(\bX_1^\top \vep - \lambda \tilde s_1)
	\]
	If we define $\hat \beta : = \{\hat \beta_1 , {\bf 0} \}$, as above, then we conclude the following
	\begin{align}
		\hat \mcS = \mcS \qquad \text{if} \qquad
		\begin{cases}
			\hat \beta_1 \overset{s}{=} \beta_{1}  & \text{and} \\
			\|X_i^\top (Y -  \bX_1 \hat \beta_1)\| < \lambda \tilde w_i & \forall \  i \notin \mcS.
		\end{cases} \label{e:prop1}
	\end{align}
In the scalar case, the first condition above is replaced with the property that $\hat \beta_i$ and $\beta_i$ point in the same direction (as this will imply that $\hat \beta_1 \overset{s}{=} \beta_{1}$).  Since we are working in the functional setting, we cannot use such a replacement.  Instead we use the slightly stronger property $\| \hat \beta_i - \beta_i \| <\| \beta_i\|$ for all $i \in \mcS$.  This says that the estimate lies in a ball centered at $\beta_i$ with radius less than $\beta_i$.  In the scalar case this will imply that the two have the same sign.
	
	Examining the second condition in \eqref{e:prop1} we have that
\begin{align*}
Y -  \bX_1 \hat \beta_1 =&  \bX_1 \beta_1  + \vep - \bX_1\beta_1 +  \bX_1 (\bX_1^\top \bX_1)^{-1}(\bX_1^\top \vep - \lambda \tilde s_1) \\
=& (\bI - \bX_1 (\bX_1^\top \bX_1)^{-1}\bX_1^\top ) \vep -  \lambda \bX_1 (\bX_1^\top \bX_1)^{-1} \tilde s_1.	
\end{align*}
Let $\bH =  (\bI - \bX_1 (\bX_1^\top \bX_1)^{-1}\bX_1^\top )$ then we get the following
\begin{align}
\hat \mcS = \mcS \qquad \text{if} \qquad
\begin{cases}
\| e_i^\top ( N^{-1} \hat \Sigma_{11}^{-1}(\bX_1^\top \vep - \lambda \tilde s_1))\| < \| \beta_i  \|  &  \forall \ i \in \mcS \\
			\|X_i^\top (\bH \vep
			-  \lambda N^{-1} \bX_1 \hat \Sigma_{11}^{-1} \tilde s_1)\| < \lambda \tilde w_i & \forall \ i \notin \mcS,
		\end{cases} \label{e:prop2}
	\end{align}
	where $e_i$ is a vector with $1$ in the $i^{th}$ coordinate and zero everywhere else.  Now define the following four events:
	\begin{align*}
		B_1 &= \left\{\frac{1}{N} \| e_i^\top \Sigma_{11}^{-1}\bX_1^\top \vep \| \geq \frac{\| \beta_i\|}{2}: 
		\ \text{for some } i \in \mcS \right\} \\
		B_2 &= \left\{\frac{\lambda}{N} \| e_i^\top \Sigma_{11}^{-1}\tilde s_1 \| \geq \frac{\| \beta_i\|}{2}: 
		\ \text{for some } i \in \mcS \right\} \\
		B_3 &= \left\{\| X_i^\top \bH \vep  \| \geq  \frac{\lambda \tilde{w_i}}{2}: 
		\ \text{for some } i \notin \mcS \right\} \\
		B_4 &= \left\{\frac{1}{N}\| X_i^\top \bX_1 \Sigma_{11}^{-1} \tilde s_1  \| \geq   \frac{ \tilde{w_i}}{2}: 
		\ \text{for some } i \notin \mcS \right\}.
	\end{align*}

We have the inclusion
	\[
	\{\hat \mcS \neq \mcS\}
	\subset B_1 \cup B_2 \cup B_3 \cup B_4.
	\]
So we can show that $P(\hat \mcS \neq \mcS) \to 0$ if we can show that $P(B_{i}) \to 0$ for $j=1,\dots,4$.   Showing this for $B_1$ and $B_3$ is very similar to the scalar case.   For $B_2$ and $B_4$ this will require more work as our $\tilde s_1$ is a bit different than in the scalar case since we can't assume that $\hat \beta_i$ and $\beta_i$ have the same sign (for $i \in \mcS$). \\

\noindent {\bf Step 1 $P(B_1) \to 0$: }  
Denote by $e_i \in \mbR^{I_0}$ a vector which is 1 in the i$^{th}$ coordinate and zero everywhere else.  Then we can express $B_{1}=\cup_{i\in S} A_{i}$ where
\begin{align*}
A_{i}&=\left\{\frac{1}{N}\|e_{i}^{T}\Sigma_{11}^{-1}X_{1}^{T}\varepsilon\|\geq \frac{\|\beta_{i}\|}{2}\right\}
 =\left\{N^{-2} \| {T_{i}}\|^{2}\geq \frac{\|\beta_{i}\|^{2}}{4}\right\},
\end{align*}
where $T_i = e_{i}^{T}\Sigma_{11}^{-1}X_{1}^{T}\varepsilon$.  Since the $\vep$ are Gaussian we have that $\frac{1}{N}T_{i}\sim \mcN(0, \frac{1}{N}e_{i}^{T}\Sigma_{11}^{-1}e_{i}C)$
Trivially we have $P(B_{1})\leq \sum P(A_{i})$.  
Applying Lemma \ref{l:barb},  we have that for any $t > 0$
\begin{equation*}
 P\left(N^{-2} \|{T_{i}}\|^{2}\geq \frac{1}{N}\left(e_{i}^{T}\Sigma_{11}^{-1}e_{i}\|C\|_{1}+2e_{i}^{T}\Sigma_{11}^{-1}e_{i}\|C\|_{2}\sqrt{t}+2e_{i}^{T}\Sigma_{11}^{-1}e_{i}\|C\|_{\infty}t\right)\right)\leq \exp\left(-t\right).
\end{equation*}
We can apply this to our problem if we can find $t$ such that
\begin{equation} \label{e:bound1}
\frac{1}{N}\left(e_{i}^{T}\Sigma_{11}^{-1}e_{i}\|C\|_{1}+2e_{i}^{T}\Sigma_{11}^{-1}e_{i}\|C\|_{2}\sqrt{t}+2e_{i}^{T}\Sigma_{11}^{-1}e_{i}\|C\|_{\infty}t \right)\leq\frac{b_{N}^{2}}{4},
\end{equation}
where recall that $b_N = \min_{i \in S} \| \beta_i\|$.  Applying Assumption \ref{a:tech}.3 we have that
	\begin{equation*}
	\frac{1}{N}e_{i}^{T}\Sigma_{11}^{-1}e_{i}\left(\|C\|_{1}+2\|C\|_{2}\sqrt{t}+2\|C\|_{\infty}t\right)\leq \frac{tK \tau_1}{N }
	\end{equation*}
where  K is some constant  depending only on $C$, the covariance operator of the errors.  So \eqref{e:bound1} holds if 
\begin{equation*}
\frac{tK \tau_1}{N}\leq\frac{b_{N}^{2}}{4}.
\end{equation*}
So we take $t=\frac{b_{N}^{2} N}{4K \tau_1}$
, which gives
\begin{equation*}
P(A_{i})=P\left(\frac{1}{N^2}\|T_{i}\|^2\geq\frac{b_{N}^{2}}{4}\right)\leq \exp\left(-\frac{b_{N}^{2}N}{4K \tau_1}\right).
\end{equation*}
And, \begin{align*}
P(B_1)&\leq\sum P\left(A_{i}\right) 
\leq I_{0} \exp\left(-\frac{b_{N}^{2}N}{4K \tau_1}\right) 
= \exp\left(-\frac{b_{N}^{2}N}{4K \tau_1}+ \log(I_0) \right) 
\to 0, 
\end{align*}
by Assumption \ref{a:tech}.1, thus, $P(B_1) \to 0$, which gives the desired result. \\

\noindent { \bf Step 2 $P(B_2) \to 0$:} 
Recall,  
\[
B_2 = \left\{\frac{\lambda}{N} \| e_i^\top \Sigma_{11}^{-1}\tilde s_1 \| \geq \frac{\| \beta_i\|}{2}: \ \text{for some } i \in \mcS \right\}, \\
\]
and $\tilde s_{1} = \{ \tilde w_{i} \| \hat \beta_{i} \|^{-1} \hat \beta_{i}, i \in \mcS\}$.   Trivially we have	
\begin{align*}
\|\tilde s_1\|^{2}
=&\sum_{i \in \mcS} \tilde w_{i}^{2}\frac{\|\hat\beta_{i}\|^{2}}{\|\hat \beta_{i}\|^{2}}
= \sum_{i \in \mcS}\tilde w_{i}^{2}
= \sum_{i \in \mcS} w_{i}^{2} + \sum_{i \in \mcS}(\tilde w_{i}^{2} - w_i^2).
\end{align*}	
Recall that $\tilde w_i = \| \tilde \beta_i\|^{-1}$, where $\tilde \beta_i$ is computed from FSL.  For a differentiable functional $f:\mathcal{H}\to \mathbb{R}$, a first order Taylor expansion is given by \\
\begin{equation*}
f(x+h)=f(x)+\langle h, f^{'}(x)\rangle+O(h^{2}).
\end{equation*}
where $f^{'}: \mathcal{H}\to \mathcal{H}$. In our case, 
\begin{equation*}
f(x)=\frac{1}{\|x\|^2},
\quad \text{and} \quad
f^{'}(x)=-\frac{2}{\|x\|^{4}} x ,~~x\neq 0.
\end{equation*}
Applying a first-order Taylor expansion for $\tilde w_{i}^2$, we obtain
		\begin{align*}
		\tilde w_{i}^2- w_{i}^2&=\frac{1}{\| \tilde \beta_{i}\|^2}-\frac{1}{\|\beta_{i}\|^2} 
		 \approx\langle \tilde \beta_{i}-\beta_{i}, -\frac{2}{\|\beta_{i}\|^{4}} \beta_{i}\rangle
		=-2 \|\beta_{i}\|^{-4}\langle \beta_{i}-\tilde\beta_{i}, \beta_{i} \rangle.
		\end{align*}
Then, by the Cauchy--Schwarz inequality,
		\begin{align*}
\frac{2}{\|\beta_{i}\|^{4}}|\langle\beta_{i}-\tilde \beta_{i},\beta_{i}\rangle|
\leq \frac{2}{\|\beta_{i}\|^{3}}\|\beta_{i}-\tilde \beta_{i}\|
		&\leq \frac{1}{b_{N}}(\sup_{i\in \mathcal{S}} \|\beta_{i}-\tilde \beta_{i}\|)\frac{1}{\|\beta_{i}\|^2}
		= \frac{1}{b_{N}}O_{p}(r_{N}^{1/2})w_{i}^2,
		\end{align*}
Therefore,
	\begin{align*}
	\left| \sum \tilde w^2_{i}- w^2_{i} \right|&\leq \frac{1}{b_{N}}O_{p}(r_{N}^{1/2})\sum w^2_{i},
	\end{align*}
and since $r_N^{1/2}/b_N \to 0$ by Assumption \ref{a:tech}.1 we have that
	\begin{align}
\|\tilde s_{1}\|^2	& \leq  \left(\sum w^2_{i}\right)(1+o_{p}(1))
\leq \frac{I_0}{b_N^2}(1+o_{p}(1)). \label{e:w2}
	\end{align}
Returning to the original objective, we combine \eqref{e:w2}, an operator inequality, and Assumptions \ref{a:tech}.2 and \ref{a:tech}.3 to obtain
\begin{align*}
\frac{\lambda \| e_i^\top \Sigma_{11}^{-1}\tilde s_1 \|}{N \|\beta_i\|}  & \leq \frac{\lambda \| e_i^\top \Sigma_{11}^{-1} \|_{op} \|\tilde s_1 \|}{N \|\beta_i\|}
{\leq} \frac{\lambda \tau_1 \|\tilde s_1 \|}{N b_{N}} 
\leq \frac{\tau_1 \lambda I_0^{1/2}}{N b_{N}^2}(1+o_p(1))  
\to 0.
\end{align*}
Since this holds uniformly in $i$ we have that $P(B_2) \to 0$, which gives the desired result.\\

\noindent {\bf Step 3 $P(B_3) \to 0$:} 
Recall,
	\[
	B_3 = \left\{\| X_i^\top \bH \vep  \| \geq  \frac{\lambda \tilde{w_i}}{2} : 
	\ \text{for some } i \notin \mcS \right\} \\
	\]
where
\begin{equation*}
\bH =  (\bI - \bX_1 (\bX_1^\top \bX_1)^{-1}\bX_1^\top ),
\end{equation*}
is the orthogonal projection.
Let $B_{3}=\cup_{i\in S^{c}} A_{i}$, where
\begin{equation*}
A_{i}=\left\{\| X_i^\top \bH \vep  \| \geq  \frac{\lambda \tilde{w_i} }{2}  \right\}.
\end{equation*}
Next define, $T_{i}=X_{i}^{T}\bH\varepsilon$.  Then $T_{i}\in \mathcal{H}$ is Gaussian, has mean zero,  and covariance $C_{T_{i}}=X_{i}^{T}\bH X_{i}C$ since $\bH$ is idempotent.
%
%
Notice that we have
\begin{align*}
\| X_i^\top \bH \vep  \| &\geq  \frac{\lambda\tilde{w_i}}{2}
\Longleftrightarrow \|\tilde \beta_{i}\|\|  X_i^\top \bH \vep  \|  \geq  \frac{\lambda}{2}.
\end{align*}
By Lemma \ref{l:fsl}, 
$
\sup_{i \in S^c} \|\tilde {\beta_{i}}\| =O_{p}(r_N^{1/2})$.  
Then $A_i$ can be expressed as  
\begin{equation*}
 A_i =  \left\{ O_{p}(r_N^{1/2}) \|  X_i^\top \bH \vep  \| \geq  \frac{\lambda}{2}\right\},
 \end{equation*}
where $O_{p(1)}$ is uniform with respect to $i$.
Then, for any $\epsilon > 0$ we can find  $K >0$ large such that
\begin{align*}
	P(B_{3})&\leq \sum P(A_{i})
	\leq \epsilon/2 + \sum_{i \in S^c} \  P \left ( \|  X_i^\top \bH \vep  \|\geq  \frac{\lambda }{2K r_{N}^{1/2}} \right).
\end{align*}
So we need only to show that we can make the second term above small as well.  Again, we will apply Lemma 1, so we want to find t, such that
\begin{equation} \label{e:b3:bound}
\left(X_{i}^{T}\bH X_{i}\|C\|_{1}+2X_{i}^{T}\bH X_{i}\|C\|_{2}\sqrt{t}+2X_{i}^{T}\bH X_{i}\|C\|_{\infty}t\right) \leq \left(\frac{\lambda}{2 K r_{N}^{1/2}} \right)^{2}.
	\end{equation}
Notice that because $\bH$ is a projection and the $X_i$ are standardized, we have that
	\begin{equation*}
	X_{i}^{T}\bH X_{i}(\|C\|_{1}+2\|C\|_{2}\sqrt{t}+2\|C\|_{\infty}t)\leq N{K_2} t,
	\end{equation*}
where  $K_2$ is a second constant depending on $C$ only. 
So \eqref{e:b3:bound} holds if 
\begin{equation*}
tK_2 N \leq \left(\frac{\lambda}{2 K r_{N}^{1/2}} \right)^{2}.
\end{equation*}
Combining constants and labeling it $K_3$, we can take t=$\frac{\lambda^2}{K_3 N r_N}$, which gives
\begin{equation*}
P \left ( \|  X_i^\top \bH \vep  \|\geq  \frac{\lambda }{2K r_{N}^{1/2}} \right)
\leq \exp\left(-\frac{\lambda^2}{K_3 N r_N}\right). 
\end{equation*}
We can then bound $P(B_i)$ as
\begin{align*}
P(B_{i})&\leq\sum P\left(A_{i}\right)
\leq I  \exp\left(-\frac{\lambda^2}{K_3 N r_N}\right)+ \epsilon/2
= \exp\left(-\frac{\lambda^2}{K_3 I_0 \log(I) } + \log(I) \right)+ \epsilon/2.
\end{align*}
By assumption \ref{a:tech}.2,  we can take $N$ large to make the first term above smaller than $\epsilon/2$ as well, which gives the desired result. \\

\noindent{\bf Step 4 $P(B_4) \to 0$: } 
Recall,
		\[
		B_4 = \left\{\frac{1}{N}\| X_i^\top \bX_1 \Sigma_{11}^{-1} \tilde s_1  \| \geq  \frac{ \tilde{w_i}}{2}: 
		\ \text{for some } i \notin \mcS \right\}.\\
		\]
From \eqref{e:w2} and Assumption \ref{a:tech}.3 we have
	\begin{align*}
\|\tilde s_{1}\| & \leq \frac{I_0^{1/2}}{b_N} (1 + o_P(1))
\qquad \text{and} \qquad \|\Sigma_{11}^{-1}\|_{op} \leq \tau_1.
	\end{align*}
From Step 3 we have
\[
\sup_{i \in S^c} \tilde w_i^{-1} = \sup_{i \in S^c}  \| \tilde \beta_i\| = O_P(r_N^{1/2}).
\]
By Assumption \ref{a:tech}.4 we have
\begin{align*}
\frac{\| X_i^\top \bX_1 \Sigma_{11}^{-1}  \|_{op} }{N  }\leq \phi.
\end{align*}
Combining these together we have
\begin{align*}
\frac{2\| X_i^\top \bX_1 \Sigma_{11}^{-1} \tilde s_1  \| }{N  \tilde w_i }  & \leq \frac{2\| X_i^\top \bX_1\Sigma_{11}^{-1}  \|_{op}  \|\tilde s_{1}\|
}{N  \tilde w_i }  \\ 
&\leq \frac{\phi r_N^{1/2} I_0^{1/2}}{b_N}O_{p}(1)    \\
&= \frac{  I_0 \log^{1/2}(I)}{ b_N N^{1/2} }O_{p}(1) 
\to 0,
\end{align*}
by Assumptions \ref{a:tech}.1 and \ref{a:tech}.2  uniformly in i, thus, $P(B_4) \to 0$, which gives the desired result.\\

\noindent {\bf Proof of Theorem \ref{t:oracle}.2: }\\ 
As we have just shown, when $\hat \mcS = \mcS$ then the AFSL estimator takes the form $\hat \beta = (\hat \beta_1, {\bf 0})$ where $\hat \beta \overset{s}{=} \beta$ and
\begin{align*}
\hat \beta_1 &=(\bX_{1}^{T}\bX_{1})^{-1}(\bX_{1}^{T}Y-\lambda\tilde s_{1})
=\beta_{1}+(\bX_{1}^{T}\bX_{1})^{-1}(\bX_{1}^{T}\varepsilon-\lambda\tilde s_{1}).
\end{align*}
The oracle estimator is given by $\hat \beta_{O} = (\hat \beta_{1O}, {\bf 0})$ where
\begin{align*}
	\hat \beta_{1O} &= (\bX_1^\top \bX_1)^{-1} (\bX_1^\top Y )
	= \beta_1 +  (\bX_1^\top \bX_1)^{-1}(\bX_1^\top \vep).
	\end{align*}
We then have, for $\epsilon > 0$
\begin{align*}
P\left\{\sqrt{N}\| \hat \beta_{O}-{\hat \beta}\| \geq \epsilon \right \}
\leq P(\hat \mcS \neq \mcS) +P\left\{\sqrt{N}\| \hat \beta_{1O}-{\hat \beta_{1}}\| \geq \epsilon \right\}.
\end{align*}
The first term can be made arbitrarily small by Theorem \ref{t:oracle}.1.  Turning to the second term
\begin{align*}
\sqrt{N}\| \hat \beta_1-\hat \beta_{1O}\| 
=\frac{\lambda}{\sqrt{N}} \| \Sigma_{11}^{-1}\tilde s_{1} \| .
\end{align*}
From \eqref{e:w2} and Assumption \ref{a:tech2}we have
\begin{align*}
\frac{\lambda}{\sqrt{N}} \| \Sigma_{11}^{-1}\tilde s_{1} \| 
\leq \frac{\lambda \tau_1 I_0^{1/2}}{\sqrt{N}b_N}O_{p}(1) 
\to 0,
\end{align*}
which gives the desired result.

\bibliographystyle{oxford3}
\bibliography{references}

\end{document}